\documentclass[12pt]{article}

\newcommand{\dfr}[2]{\dfrac{#1}{#2}}
\newcommand{\cd}{\cdot}
\newcommand{\cds}{\cdots}
\newcommand{\dsum}{\displaystyle \sum}

\newcommand{\ol}[1]{\overline{#1}}
\newcommand{\ul}[1]{\underline{#1}}

\renewcommand{\l}{\left}
\renewcommand{\r}{\right}

\newcommand{\vsv}{\vspace{5mm}}
\newcommand{\vsb}{\vspace{2mm}}
\newcommand{\q}{\quad}
\newcommand{\qq}{\qquad}

\newcommand{\la}{\langle}
\newcommand{\ra}{\rangle}

\newcommand{\abs}[1]{\lvert{#1}\rvert}

\newcommand{\Z}{\mathbb{Z}}
\newcommand{\C}{\mathbb{C}}
\newcommand{\R}{\mathbb{R}}
\newcommand{\N}{\mathbb{N}}

\newcommand{\M}{\mathbb{M}}

\newcommand{\vir}{\mathrm{Vir}}
\newcommand{\aut}{\mathrm{Aut}}

\renewcommand{\hom}{\mathrm{Hom}}

\newcommand{\id}{\mathrm{id}}

\newcommand{\om}{\omega}

\newcommand{\al}{\alpha}
\newcommand{\ts}{\tilde{s}}

\newcommand{\pii}{\pi \sqrt{-1}\, }

\newcommand{\Span}{\mathrm{span}}

\newcommand{\w}{\omega}
\newcommand{\vacuum}{\mathbbm{1}}
\newcommand{\vac}{\vacuum}
\newcommand{\irr}{\mathrm{Irr}}

\newcommand{\Fi}{\mathrm{Fi}_{24}}
\newcommand{\BN}{\mathcal{G}^\natural}

\makeatletter
\@addtoreset{equation}{section}

\makeatother

\theoremstyle{plain}
\newtheorem{mthm}{Main Theorem}

\theoremstyle{definition}

\title{McKay's $E_6$ observation on the largest Fischer group}
\author{
  Gerald H\"{o}hn\footnote{Partially supported by Kansas NSF EPSCoR grant NSF32239KAN32240.}
  \vsb\\
  {\small \it Department of Mathematics, Kansas State University, }
  \\
  {\small \it 138 Cardwell Hall, Manhattan, KS 66506-2602, USA}
  \\
  {\small e-mail: {\tt gerald@math.ksu.edu}}
  \vsv\\
  Ching Hung Lam\footnote{Partially supported by NSC grant
  97-2115-M-006-015-MY3  and National Center for Theoretical Sciences, Taiwan.}%
  \vsb\\
  {\small \it Institute of Mathematics, Academia Sinica, Taipei, Taiwan 11529}\\
  {\small e-mail: {\tt chlam@math.sinica.edu.tw}}
  \vsv\\
  Hiroshi Yamauchi\footnote{Partially supported by JSPS Grant-in-Aid
  for Young Scientists (Start-up) No.\ 19840025 and (B) No.\ 21740011.}
  \vsb\\
  {\small \it Department of Mathematics,
  Tokyo Woman's Christian University}\\
  {\small \it 2-6-1 Zempukuji, Suginami-ku, Tokyo 167-8585, Japan}\\
  {\small e-mail: {\tt yamauchi@lab.twcu.ac.jp}}
  \vsv
}
\date{}

\newcommand{\sfr}[2]{\leavevmode\kern-.05em
  \raise.5ex\hbox{\the\scriptfont0 #1}\kern-.1em
  /\kern-.15em\lower.25ex\hbox{\the\scriptfont0 #2}\kern.02em}

\DeclareMathOperator*{\tensor}{\otimes}
\DeclareMathOperator*{\fusion}{\times}

\newcommand{\pf}{\noindent {\bf Proof:}\q }
\newcommand{\stab}{\mathrm{Stab}}

\newcommand{\B}{\mathbb{B}}
\newcommand{\com}{\mathrm{Com}}

\newcommand{\VF}{V\!\! F^\natural}
\newcommand{\W}{\EuScript{W}}
\newcommand{\EL}{\tilde{L}}

\newcommand{\VFnX}[1]{V_{F(#1)}}

\newcommand{\UFnX}[1]{U_{F(#1)}}
\newcommand{\GFnX}[1]{\mathcal{G}_{F(#1)}}

\newcommand{\tv}{\tilde{v}}

\begin{document}

\baselineskip 6mm

\maketitle

\thispagestyle{empty}

\begin{abstract}
In this paper, we study McKay's $E_6$-observation on the largest Fischer
$3$-trans\-position group $\Fi$.
We investigate a vertex operator algebra $\VF$ of central charge $23\frac{1}{5}$
on which the Fischer group $\Fi$ naturally acts.
We show that there is a natural correspondence between dihedral subgroups of
$\Fi$ and certain vertex operator subalgebras constructed by the nodes of the
affine $E_6$ diagram by investigating so called derived Virasoro vectors of
central charge~$6/7$.
This allows us to reinterpret McKay's $E_6$-observation via the theory of
vertex operator algebras.

It is also shown that the product of two non-commuting Miyamoto involutions
of $\sigma$-type associated to derived $c=6/7$ Virasoro vectors is an element
of order~$3$, under certain general hypotheses on the vertex operator algebra.
For the case of $\VF$, we identify these involutions
with the $3$-transpositions of the Fischer group~$\Fi$.

\end{abstract}

\renewcommand{\arraystretch}{1.5}

\newpage

\vbox{\tableofcontents}

\section{Introduction}

This article is a continuation of our previous work~\cite{LM,LYY1,HLY} to give
a vertex operator algebra (VOA) theoretical interpretation of McKay's intriguing
observations that relate the Monster, the Baby Monster and the largest Fischer
$3$-transposition group to the affine $E_8$, $E_7$ and $E_6$ Dynkin diagrams.
In this article, we will study the $E_6$-observation recalled below.
Our approach here is similar to \cite{HLY}, in which the $E_7$-observation
is studied, but other vertex operator algebras are involved and many
technical details are different.


\paragraph{The largest Fischer group.}

The largest Fischer group $\Fi$ was discovered by B.~Fischer~\cite{Fischer-not}
as a group of order $2^{22} \cdot 3^{16} \cdot 5^2 \cdot 7^3 \cdot 11 \cdot 13
\cdot 17 \cdot 23 \cdot 29$ containing a conjugacy class of $306,936$
involutions, which satisfy the {\it $3$-transposition property,\/} i.e.,\ any
non-commuting pair has product of order~$3$.
The group $\Fi$ contains as subgroup of index~$2$ the derived group $\Fi'$
which is the third largest of the 26 sporadic groups.
An extension $3.\Fi$ of the Fischer group $\Fi$ by a cyclic group of order $3$ is
the normalizer of a 3A-element of the Monster, the largest of the sporadic groups.
In fact, one can construct $\Fi$ from the Monster \cite{G} and derive its
$3$-transposition property.
One of the main motivation of this article is to study the Fischer group $\Fi$
and to understand the $E_6$ case of McKay's observations by using the theory of
vertex operator algebra.

There exists a class of vertex operator algebras which is closely related to
$3$-transposition groups~\cite{Gr3,KM,M1,Ma2}.
Let $V$ be vertex operator algebra which has a simple $c=1/2$ Virasoro vector
$e$ of {\it $\sigma$-type}, that is $V= V_e[0]\oplus V_e[\sfr{1}{2}]$
where $V_e[h]$ denotes the sum of irreducible $\vir(e)$-submodules of $V$
isomorphic to the irreducible highest weight representation $L(\sfr{1}{2},h)$ of
the $c=1/2$ Virasoro algebra with highest weight $h$ (cf.~Section~\ref{sec2}).
Then one can define an involutive automorphism
\[
\sigma_e=
\begin{cases}
1 & \text{ on } V_e[0], \\
-1 &\text{ on } V_e[\sfr{1}{2}]\\
\end{cases}
\]
usually called a $\sigma$-involution if $\sigma_e\not=\id_V$~\cite{M1}.
It was shown by Miyamoto \cite{M1} that if the weight
one subspace of a VOA is trivial, then a collection of involutions associated to
$c=1/2$ Virasoro vectors of $\sigma$-type generates a $3$-transposition group.
Many interesting examples of $3$-transposition groups obtained by $\sigma$-type $c=1/2$
Virasoro vectors have been studied in \cite{Gr3,KM} and the
complete classification is established in \cite{Ma2}. According to \cite{Ma2}, all
$3$-transposition groups realized by $\sigma$-type $c=1/2$  Virasoro vectors are
so-called  {\it symplectic type\/} (cf.\ \cite{CH}), and as a result, the Fischer
$3$-transposition group cannot be obtained by $c=1/2$ Virasoro vectors of
$\sigma$-type.

Another result on $c=1/2$ Virasoro vectors was obtained in~\cite{S} where
it was shown that the so-called $\tau$-involutions (cf.\ Theorem~\ref{thm:2.2})
associated to such vectors generate a $6$-transposition group
provided the weight one subspace of the VOA is trivial.

\smallskip

In this paper, we introduce a new idea to obtain $3$-transposition groups as
automorphism groups of vertex operator algebras. We use the so-called
3A-algebra for the Monster (see \cite{M3,LYY2,SY} and Section \ref{sec:3})  and
consider derived $c=6/7$ Virasoro vectors to define involutive automorphisms of
vertex operator algebras. We will show that a collection of involutions associated
to derived $c=6/7$ Virasoro vectors generates a $3$-transposition group. The
advantage of our method is that we can realize the largest Fischer
$3$-transposition group as an automorphism subgroup of a special vertex
operator algebra $\VF$ explained below. This result enables us to study the
Fischer $3$-transposition groups via the theory of vertex operator algebras.

\paragraph{The Fischer group VOA $\VF$.}

We will investigate a certain vertex operator algebra $\VF$ of central charge
$23\frac{1}{5}$ on which the Fischer group $\Fi$ naturally acts.

\smallskip

The Monster is the automorphism group  of the Moonshine vertex operator
algebra $V^\natural$ (cf.~\cite{FLM,B}).
Let $g$ be a 3A-element of the Monster $\M$.
Then the normalizer $N_\M(\la g\ra)$ is isomorphic to $3.\Fi$ and acts on $V^\natural$.
A character theoretical consideration in \cite{C,MN} indicates that
the centralizer $C_{\M}(g)\simeq 3.\Fi'$ fixes a unique $c=4/5$ Virasoro vector
in $V^\natural$.
We will show in Theorem \ref{thm:5.16} that $C_\M(g)$ actually fixes a unique
$c=4/5$ extended Virasoro vertex operator algebra
$\W\simeq \W(\sfr{4}{5})=L(\sfr{4}5,0)\oplus L(\sfr{4}5,3)$ in $V^\natural$.
Let
\[ \VF =\com_{V^\natural}(\W), \]
where $\com_V(U)$ denotes the commutant subalgebra of $U$ in $V$
(see~\eqref{eq:3.1} and \eqref{eq:5.12} for the precise definition),
and we call $\VF$ the {\it Fischer group VOA}.
A simple observation shows that $N_\M(\la g\ra)$ acts naturally on $\VF
=\com_{V^\natural}(\W)$.
In fact, we will show that the Fischer group $\Fi$ can be realized as a
subgroup of $\aut(\VF)$.

\begin{mthm}[Theorem \ref{thm:5.19}]\label{mthm:1}
  The automorphism group $\aut(\VF)$ of $\VF$ contains $\Fi$ as a subgroup.
  Moreover, let $\EuScript{X}$ be the full-subalgebra of $\VF$ generated by
  its weight $2$ subspace.
  Then $\aut(\EuScript{X})\simeq \mathrm{Fi}_{24}$.
\end{mthm}

Although we only show that $\Fi$ equals the automorphism group of the Griess
algebra of $\VF$, we expect that $\aut(\VF)$ is exactly $\Fi$ and therefore
$\VF$ would provide a VOA model for studying the Fischer group $\Fi$.

Since $N_{\M}(\la g\ra )\simeq 3.\Fi$ for any 3A-element $g$ of the Monster, the study of the $3$-transpositions in $\Fi$
leads to the study of dihedral subalgebras~$U_{3A}$ of type 3A in $V^\natural$
(cf.\ \cite{M3,LYY2,SY}).
The 3A-algebra $U_{3A}$ contains a unique extended $c=4/5$ Virasoro sub-VOA
$\W=\W(\sfr{4}{5})$ and the corresponding commutant subalgebra
$\com_{U_{3A}}(\W)$ in $U_{3A}$ (cf.~\cite{SY}) gives rise to certain $c=6/7$ Virasoro vectors
in $\VF$, which we call {\it derived\/} Virasoro vectors
(cf.~Definition~\ref{df:3.13}).
Note that $\W \subset U_{3A}\subset V^\natural$ implies $\VF=
\com_{V^\natural}(\W)\supset \com_{U_{3A}}(\W)$.
It is interesting that a natural construction in~\cite{DLMN} suggests
\hbox{$c=6/7$} Virasoro vectors in the case of $E_6$; cf.~the discussion below.
Motivated by the above observation, we first study subalgebras $U_{3A}$ of
any vertex operator algebra $V$ and we show that
one can canonically associate involutive automorphisms to derived~$c=6/7$
Virasoro vectors in $\com_V(\W)$, which we call {\it Miyamoto involutions}
(see Lemma \ref{lem:2.5} and Eq.~\eqref{eq:2.4}).
We will show that the collection of Miyamoto involutions associated to
derived $c=6/7$ Virasoro vectors generates a $3$-transposition group.

We say that a VOA $W$ over $\R$ is {\it compact} if $W$ has a positive definite
invariant bilinear form.
A real sub-VOA $W$  is said to be a {\it compact real form} of a VOA $V$ over $\C$
if $W$ is compact and $V\simeq \C\otimes_\R W$.

\begin{mthm}[Theorem~\ref{prop:3.16}]\label{mthm:2}
  Let $V=\bigoplus_{\geq 0} V_n$ be a VOA.
  Suppose that $\dim V_0=1$, $V_1=0$ and $V$ has a compact real form $V_\R$
  and every simple $c=1/2$ Virasoro vector of $V$ is in $V_\R$.
  Then the Miyamoto involutions associated to derived $c=6/7$ Virasoro vectors
  in the commutant subalgebra $\com_V(\W(\sfr{4}{5}))$ satisfy
  a $3$-transposition property.
\end{mthm}

The Moonshine vertex operator algebra satisfies the assumption of the theorem
above and we recover in a general fashion the $3$-transposition property of
the largest Fischer group via the commutant subalgebra
$\VF =\com_{V^\natural}(\W(\sfr{4}{5}))$ (see Corollary \ref{cor:5.25}).
Indicated by the fact that 3-transpositions are induced by $c=6/7$ derived
Virasoro vectors, we will also prove the following one-to-one correspondence.



\begin{mthm}[Theorem \ref{thm:5.25}]\label{mthm:3}
  There exists a one-to-one correspondence between $3$-transpositions of
  the Fischer group and derived $c=6/7$ Virasoro vectors in $\VF$
  via Miyamoto involutions.
\end{mthm}

This theorem provides a link for studying the $3$-transpositions of $\Fi$
by using the VOA $\VF$ and it is possible to relate McKay's $E_6$-observation
of $\Fi$ to the theory of vertex operator algebra as discussed below.


%
%


\paragraph{McKay's observation.}

We are interested in the $E_6$ case of the observation of
McKay which relates the Monster group and some sporadic groups involved in
the Monster to the affine Dynkin diagrams of types $E_6$, $E_7$ and
$E_8$ \cite{Mc}.

McKay's observation in the case of the largest Fischer group says that the
orders of the products of any two  $3$-transpositions of $\Fi$ belongs to one
of the conjugacy classes 1A, 2A or 3A of $\Fi$ such that these conjugacy
classes coincide with the numerical labels of the nodes in an affine $E_6$
Dynkin diagram and there is a correspondence as follows:
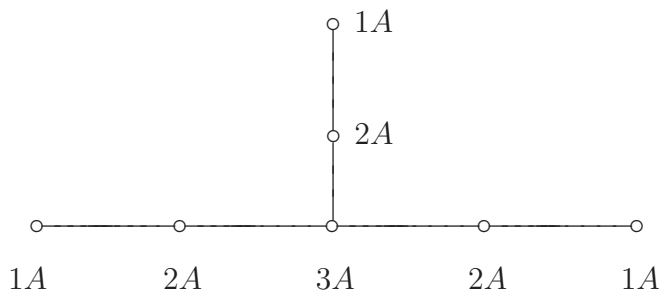
\begin{equation*}\label{eq:mckayE6+}
\begin{array}{l}
  \hspace{120pt}\circ \hspace{5pt}  1A
  \vspace{-13.5pt}\\
  \hspace{121.3pt}| \vspace{-13pt}\\
  \hspace{121.3pt}| \vspace{-13pt}\\
  \hspace{121.3pt}| \vspace{-13pt}\\
  \hspace{121.3pt}| \vspace{-13.5pt}\\
  \hspace{120.1pt}\circ \hspace{5pt}  2A
  \vspace{-13.5pt}\\
  \hspace{121.3pt}| \vspace{-13pt}\\
  \hspace{121.3pt}| \vspace{-13pt}\\
  \hspace{121.3pt}| \vspace{-13.5pt}\\
  \hspace{7.8pt}
  \circ\hspace{-5pt}-\hspace{-7pt}-\hspace{-7pt}-\hspace{-7pt}-
  \hspace{-7pt}-\hspace{-7pt}\hspace{-7pt}-\hspace{-7pt}-\hspace{-7pt}-\hspace{-7pt}-
  \hspace{-7pt}-\hspace{-6pt}-\hspace{-5pt}
  \circ\hspace{-5pt}-\hspace{-7pt}-\hspace{-7pt}-\hspace{-7pt}-
  \hspace{-7pt}-\hspace{-7pt}\hspace{-7pt}-\hspace{-7pt}-\hspace{-6pt}-\hspace{-6pt}-
  \hspace{-6pt}-\hspace{-6pt}-\hspace{-5pt}
  \circ\hspace{-5pt}-\hspace{-7pt}-\hspace{-7pt}-\hspace{-7pt}-
  \hspace{-7pt}-\hspace{-7pt}\hspace{-7pt}-\hspace{-7pt}-\hspace{-6pt}-\hspace{-6pt}-
  \hspace{-6pt}-\hspace{-6pt}-\hspace{-5pt}
  \circ\hspace{-5pt}-\hspace{-7pt}-\hspace{-7pt}-\hspace{-7pt}-
  \hspace{-7pt}-\hspace{-7pt}\hspace{-7pt}-\hspace{-7pt}-\hspace{-6pt}-\hspace{-6pt}-
  \hspace{-6pt}-\hspace{-6pt}-\hspace{-5pt}
  \circ
  \vspace{-1pt}\\
  1A\hspace{44pt}  2A \hspace{43  pt} 3A\hspace{43pt} 2A\hspace{43pt}
  1A
  \end{array}
\end{equation*}
This correspondence is not one-to-one but only up to
diagram automorphisms.

\smallskip

For the understanding of the $E_8$-case \cite{LM,LYY1,LYY2}, the main foothold
is the one-to-one correspondence between 2A-involutions of the Monster and
simple $c=1/2$ Virasoro vectors in $V^\natural$ by which one can translate
McKay's $E_8$-observation into a purely vertex operator algebra theoretical
problem. For the $E_6$ case, we have a nice correspondence as in Main Theorem
\ref{mthm:3} and we can also translate the $E_6$-observation into a problem of
vertex operator algebras. Based on this correspondence and making use of the
3A-algebra $U_{3A}$, which is generated by two $c=1/2$ Virasoro vectors, we
show that there is a natural connection between dihedral subgroups of $\Fi$ and
certain sub-VOAs constructed by the nodes of the affine~$E_6$ diagram, which
gives some context of 
McKay's observation in terms of vertex operator
algebras.

More precisely, let $S$ be a simple root lattice with a simply laced root
system $\Phi(S)$.
We scale $S$ such that the roots have squared length~$2$.
Let  $V_{\sqrt{2}S}$ be the lattice VOA associated to $\sqrt{2}S$.
Here and further we use the standard notation for lattice VOAs as in \cite{FLM}.
In \cite{DLMN} Dong et al.\ constructed a Virasoro vector of $V_{\sqrt{2}S}$ of
the form
\begin{equation}\label{eq:1.2}
  \tilde{\w}_S
  := \dfr{1}{2h(h+2)}\dsum_{\alpha \in \Phi(S)} \alpha(-1)^2 \vac
  +\dfr{1}{h+2}\dsum_{\alpha\in \Phi(S)} e^{\sqrt{2}\alpha},
\end{equation}
where $h$ denotes the Coxeter number of $S$.
Recall that the central charge of $\tilde{\w}_S$ is 6/7 if $S=E_6$ \cite{DLMN}.
By the expression, it is clear that $\tilde{\w}_S$ is invariant under
the Weyl group of $\Phi(S)$.

Our approach to McKay's observation is to find suitable pairs of derived
$c=6/7$ Virasoro vectors in $\VF$ which inherits the $E_6$ structure in
McKay's $E_6$-diagram.
By using similar ideas as in \cite{LYY1, LYY2}, we construct a certain sub-VOA
$\UFnX{nX}$ of the lattice VOA $V_{\sqrt{2}E_6}$ associated to each node $nX$
of the affine $E_6$ diagram (cf.\ Section~\ref{sec:4}).
Utilizing an embedding of the $E_6$ lattice into the $E_8$ lattice,
we show that $\UFnX{nX}$ is contained in the VOA $\VF$ purely by their
VOA structures (cf.~Theorem \ref{thm:5.30} and Appendix A).
These sub-VOAs contain pairs of derived $c=6/7$ Virasoro vectors such that
the corresponding Miyamoto involutions generate a dihedral group of type
$nX$ in $\Fi\subset \aut(\VF)$.
Then, using the identification of $\Fi$ as a subgroup of $\aut(\VF)$,
we obtain another main result with the help of the  Atlas~\cite{ATLAS}.

\begin{mthm}[Theorem \ref{thm:5.30}]\label{mthm:4}
  For any of the cases $nX=1A$, $2A$ or $3A$, the VOA $\UFnX{nX}$
  can be embedded into $\VF$.
  Moreover, $\sigma_{\tv}\sigma_{\tv'}$ belongs to the
  conjugacy class $nX$ of $\mathrm{Fi}_{24}=\aut(\EuScript{X})$,
  where $\tv$ and $\tv'$ are defined as in \eqref{tv}.
\end{mthm}

In this way, our embeddings of $\UFnX{nX}$ into $\VF$ encode the $E_6$
structure into $\VF$ which are compatible with the original McKay observations.

\medskip


\paragraph{The organization of the paper.}

The organization of this article is as follows:
In Section 2, we review basic properties about Virasoro VOAs and
Virasoro vectors.

In Section \ref{sec:4}, we recall the definition of commutant sub-VOAs and
define certain commutant subalgebras
associated to the root lattice of type  $E_6$ using the method
described in \cite{LYY1,LYY2}.

In Section \ref{sec:3}, we study a vertex operator algebra $U_{3A}$, which we
call the 3A-algebra for the Monster, and prove a $3$-transposition property for
Miyamoto involutions associated to derived $c=6/7$ Virasoro vectors on
commutant subalgebras of VOAs containing $U_{3A}$.

In Section \ref{sec:5}, the commutant subalgebra $\VF$ of $\W(\sfr{4}{5})$ in
$V^\natural$ is studied.
The sub-VOA $\VF$ has a natural faithful action of the Fischer group $\Fi$.
We expect that the full automorphism group of $\VF$ is $\Fi$ but we cannot give
a proof of this prospect. Instead, we will show a partial result that
the full automorphism group of the sub-VOA $\EuScript{X}$ generated by the
weight two subspace of $\VF$ is isomorphic to $\Fi$. We also establish a
one-to-one correspondence between 2C-involutions of $\Fi$ and derived
$c=6/7$ Virasoro vectors in $\VF$.

Finally,  we discuss the embeddings of the
commutant subalgebras constructed in Section \ref{sec:4} into
$\VF$ in Section~\ref{sec:5.2.3}.
We show that the $c=6/7$ Virasoro vectors defined in Section~\ref{sec:4.1}
can be embedded into $\VF$.
Moreover, we verify that the product of the corresponding $\sigma$-involutions
belongs to the conjugacy class associated to the node.


\paragraph{Acknowledgment.}
Part of the work was done when the first and third named authors
were visiting the National Center for Theoretical Sciences, Taiwan
in June 2008 and when the third named author was visiting
Kansas State University in April 2006.
They gratefully acknowledge the hospitalities there received.
The first named author likes to thank Simon Norton and Alexander Ivanov for providing
him with an explicit description of the Griess algebra of $V_{F(3A)}$.
The third named author also thanks Hiroki Shimakura for valuable comments
on automorphism groups of vertex operator algebras.
All of the authors thank ICMS for the opportunity to have a discussion on
this work during the conference in Edinburgh, September 2009.
Finally, the authors thank the referee for the careful reading
of the manuscript and his/her valuable comments.


\paragraph{Notation and Terminology.}
In this article, $\N$, $\Z$, $\R$ and $\C$ denote the set of non-negative
integers, integers, real and complex numbers, respectively. We denote the ring
$\Z/p\Z$ by $\Z_p$ with a positive integer $p$ and often identify the integers
$0,1,\dots, p-1$ with their images in $\Z_p$.

Every vertex operator algebra (VOA for short) is defined over the field $\C$ of
complex numbers unless otherwise stated. A VOA $V$ is called {\it of CFT-type}
if it is non-negatively graded $V=\bigoplus_{n\geq 0}V_n$ with $V_0=\C \vac$.
For a VOA structure $(V,Y(\cd,z),\vac,\w)$ on $V$, the vector $\w$ is called
the {\it conformal vector} of $V$ \footnote{The conformal vector of $V$ is often
called the Virasoro element of $V$, e.g.~\cite{FLM}}.
For simplicity, we often use $(V,\w)$ to
denote the structure $(V,Y(\cd,z),\vac,\w)$. The vertex operator $Y(a,z)$ of
$a\in V$ is expanded  as $Y(a,z)=\sum_{n\in \Z}a_{(n)} z^{-n-1}$.

An element $e\in V$ is referred to as a
{\it Virasoro vector of  central charge\/} $c_e\in \C$ if
$e\in V_2$ and it satisfies $e_{(1)}e=2e$ and $e_{(3)}e=(c_e/2)\cd \vac$.
It is well-known that the associated modes $L^e(n):=e_{(n+1)}$,
$n\in \Z$, generate a representation of the Virasoro algebra
on $V$ (cf.\ \cite{M1}), i.e., they satisfy the commutator relation
$$
  [L^e(m),L^e(n)]
  = (m-n)L^e(m+n)+\delta_{m+n,0}\dfr{m^3-m}{12}c_e.
$$
Therefore, a Virasoro vector together with the vacuum vector generates a
Virasoro VOA inside $V$. We will denote this subalgebra by $\vir(e)$.

In this paper, we define a sub-VOA of $V$ to be a pair $(U,e)$ consisting
of a subalgebra $U$ containing the vacuum element $\vac$ and
a conformal vector $e$ for $U$ such that $(U,e)$ inherits the grading of $V$,
that is, $U=\bigoplus_{n\geq 0} U_n$ with $U_n=V_n\cap U$,
but $e$ may not be the conformal vector of $V$.
In the case that $e$ is also the conformal vector of $V$,
we will call the sub-VOA $(U,e)$ a {\it full\/} sub-VOA.

For a positive definite even lattice $L$, we will denote the
lattice VOA associated to $L$ by $V_L$ (cf.\ \cite{FLM}).
We adopt the standard notation for $V_L$ as in \cite{FLM}.
In particular, $V_L^+$ denotes the fixed point subalgebra of
$V_L$ under a lift of the $(-1)$-isometry on $L$.
The letter $\Lambda$ always denotes the Leech lattice,
the unique even unimodular lattice of rank 24 without roots.

Given a group $G$ of automorphisms of $V$, we denote by $V^G$ the fixed point
subalgebra of $G$ in $V$. The subalgebra $V^G$ is called the {\it $G$-orbifold}
of $V$ in the literature.
For a $V$-module $(M,Y_M(\cd,z))$ and $\sigma\in \aut(V)$,
we set ${}^\sigma Y_M(a,z):= Y_M(\sigma^{-1} a,z)$ for $a\in V$.
Then the {\it $\sigma$-conjugated module $\sigma\circ M$} of $M$ is defined
to be the module structure $(M,{}^\sigma Y_M(\cd,z))$.


\section{Virasoro vertex operator algebras and their extensions}\label{sec2}

For complex numbers $c$ and $h$, we denote by $L(c,h)$ the
irreducible highest weight representation of the Virasoro algebra
with central charge $c$ and highest weight $h$.
It is shown in \cite{FZ} that $L(c,0)$ has a natural structure of a
simple VOA.

\subsection{Unitary Virasoro vertex operator algebras}\label{sec:2.1}

Let
\begin{equation}\label{eq:2.1}
\begin{array}{rl}
  c_m &:= 1-\dfr{6}{(m+2)(m+3)},\qq m=1,2,\dots ,
  \vsb\\
  h_{r,s}^{(m)} &:= \dfr{\{ r(m+3)-s(m+2)\}^2-1}{4(m+2)(m+3)},\q
  1\leq s\leq r\leq m+1.
\end{array}
\end{equation}
It is shown in \cite{W} that $L(c_m,0)$ is rational and
$L(c_m,h_{r,s}^{(m)})$, $1\leq s\leq r\leq m+1$, provide all
irreducible $L(c_m,0)$-modules (see also \cite{DMZ}).
This is the so-called unitary series of the Virasoro VOAs.
The fusion rules among $L(c_m,0)$-modules are computed in
\cite{W} and given by
\begin{equation}\label{eq:2.2}
  L(c_m,h^{(m)}_{r_1,s_1})\fusion L(c_m,h^{(m)}_{r_2,s_2})
  = \dsum_{\scriptstyle i\in I \atop \scriptstyle j\in J}
    L(c_m,h^{(m)}_{\abs{r_1-r_2}+2i-1,\abs{s_1-s_2}+2j-1}),
\end{equation}
where
$$
\begin{array}{l}
  I=\{ 1\,,2,\,\dots,\,\min \{ r_1,\,r_2,\,m+2-r_1,\,m+2-r_2\}\} ,
  \vsb\\
  J=\{ 1,\,2,\,\dots,\,\min \{ s_1,\,s_2,\,m+3-s_1,\,m+3-s_2\}\} .
\end{array}
$$

\begin{df}\label{df:2.1}
  A Virasoro vector $e$ with central charge $c$ is called {\it simple\/}
  if $\vir(e)\simeq L(c,0)$.
  A simple $c=1/2$ Virasoro vector is called an {\it Ising vector.}
\end{df}

\medskip

The fusion rules among $L(c_m,0)$-modules have a canonical
$\Z_2$-symmetry and this symmetry gives rise to an involutive vertex operator algebra
automorphism which is known as Miyamoto involution.

\begin{thm}[\cite{M1}]\label{thm:2.2}
  Let $V$ be a VOA and let $e\in V$ be a simple Virasoro vector
  with central charge $c_m$.
  Denote by $V_e[h^{(m)}_{r,s}]$ the sum of irreducible
  $\vir(e)=L(c_m,0)$-submodules of $V$ isomorphic to $L(c_m,h^{(m)}_{r,s})$,
  $1\leq s\leq r\leq m+1$.
  Then the linear map
  \[
    \tau_e
    :=
    \begin{cases}
      (-1)^{r+1}
      & \text{on }\ V_e[h^{(m)}_{r,s}]\q \text{if } m \ \text{is even},
      \vsb\\
      (-1)^{s+1}
      & \text{on }\ V_e[h^{(m)}_{r,s}]\q \text{if } m \ \text{is odd},
    \end{cases}
  \]
  defines an automorphism of $V$ called the $\tau$-involution associated to $e$.
\end{thm}

Next we introduce a notion of $\sigma$-type Virasoro vectors.

\begin{df}\label{df:2.3}
For $m=1$, $2$, $\dots$, let
\[
  B^{(m)}:=
  \begin{cases}
    ~\l\{ h^{(m)}_{1,s}\mid 1 \leq s\leq m+2\r\} & \text{ if $m$ is even,}
    \vsb\\
    ~\l\{ h^{(m)}_{r,1}\mid 1 \leq r\leq m+1\r\} & \text{ if $m$ is odd.}\\
  \end{cases}
\]
A simple Virasoro vector $e\in V $ with central charge $c_m$ is said to be
of {\it $\sigma$-type} on $V$ if $V_e[h]=0$ for all $h\notin B^{(m)}$.
\end{df}

By Eq.~(\ref{eq:2.2}), the fusion rules among irreducible modules $L(c_m,h)$,
$h\in B^{(m)}$, are relatively simple.
Moreover, $B^{(m)}$ possesses a natural $\Z_2$-symmetry as follows.

\begin{lem}\label{lem:2.5}
  Let $e\in V$ be a simple $c=c_m$ Virasoro vector of $\sigma$-type.
  Then one has the isotypical decomposition
  $$
    V=\bigoplus_{h\in B^{(m)}} V_e[h]
  $$
  and the linear map $\sigma_e$ given by
  \begin{equation}\label{eq:2.3}
    \sigma_e   :=
    \begin{cases}
      (-1)^{s+1}
      & \text{on }\ V_e[h^{(m)}_{1,s}]\q \text{if } m \ \text{is even},
      \vsb\\
      (-1)^{r+1}
      & \text{on }\ V_e[h^{(m)}_{r,1}]\q \text{if } m \ \text{is odd}.
    \end{cases}
  \end{equation}
  is an automorphism of $V$.
\end{lem}
We will call the map $\sigma_e$ above the {\it Miyamoto involution (of
$\sigma$-type)\/} associated to a simple $c=c_m$ Virasoro vector $e$
of $\sigma$-type.

\begin{rem}\label{rem:2.4}
  (1) By Eq.~\eqref{eq:2.2}, $B^{(m)}$ is closed under the fusion product,
  i.e., if $h$, $h'\in B^{(m)}$ then $L(c_m,h)\fusion L(c_m,h')$ is a
  sum of irreducible modules with highest weights in~$B^{(m)}$.
  Therefore, the subspace $W=\bigoplus_{h\in B^{(m)}} V_e[h]$ forms a sub-VOA of $V$.
  Note that $e$ is of $\sigma$-type on $W$ and one can always define $\sigma_e$ as
  an automorphism of $W$.
  \\
  (2) By definition, it is clear that $\tau_e$ acts trivially on $V_e[h]$ for
  all $h\in B^{(m)}$. However, the fixed point sub-VOA of $\tau_e$ on $V$ is usually bigger
  than $W=\bigoplus_{h\in B^{(m)}} V_e[h]$.
\end{rem}

In this article, we will mainly consider the case $c_4=6/7$. In this case,
$B^{(4)}=\{0,\,  5,\, \sfr{1}{7},\, \sfr{5}{7},\, \sfr{12}{7},\, \sfr{22}{7}\}$ and
a simple $c=6/7$ Virasoro vector $e\in V$ is of $\sigma$-type on $V$ if
  $V_e[h]=0$ for $h\ne 0$, $5$, $\sfr{1}{7}$, $\sfr{5}{7}$, $\sfr{12}{7}$,
  $\sfr{22}{7}$. The corresponding $\sigma$-involution is given by
\begin{equation}\label{eq:2.4}
  \sigma_e:=
  \begin{cases}
    \, \ 1 & \mathrm{on}~~ V_e[0]\oplus V_e[\sfr{5}{7}]
    \oplus V_e[\sfr{22}{7}],
    \vsb\\
    \, -1 & \mathrm{on}~~ V_e[5]\oplus V_e[\sfr{12}{7}]\oplus
    V_e[\sfr{1}{7}].
  \end{cases}
\end{equation}

\smallskip

We also
need the following result:

\begin{lem}\label{lem:2.8}
  Let $V$ be a VOA with grading $V=\bigoplus_{n\geq 0}V_n$,
  $V_0=\C \vac$ and $V_1=0$, and let $u\in V$ be a Virasoro vector
  such that $\vir(u)\simeq L(c_m,0)$.
  Then the zero mode $o(u)=u_{(1)}$ acts on the Griess algebra
  of $V$ semisimply with possible eigenvalues $2$ and
  $h_{r,s}^{(m)}$, $1\leq s\leq r\leq m+1$.
  Moreover, if $h_{r,s}^{(m)}\ne 2$ for $1\leq s\leq r\leq m+1$
  then the eigenspace for the eigenvalue $2$ is one-dimensional,
  namely, it is spanned by the Virasoro vector $u$.
\end{lem}

\pf
See Lemma 2.6 of \cite{HLY}.
\qed


\subsection{Extended Virasoro vertex operator algebras}\label{sec:2.2}
Among $L(c_m,0)$-modules, only $L(c_m,0)$ and
$L(c_m,h_{m+1,1}^{(m)})$ are simple currents, and it is shown in
\cite{LLY} that $L(c_m,0)\oplus L(c_m,h^{(m)}_{m+1,1})$ forms a simple
current extension of $L(c_m,0)$.
Note that $h_{m+1,1}^{(m)}(= h_{1,m+2}^{(m)})=m(m+1)/4$ is an integer if
$m\equiv 0,\,3 \pmod{4}$ and a half-integer if $m\equiv 1,\,2 \pmod{4}$.

\begin{thm}[\cite{LLY}]\label{thm:2.7}
  The $\Z_2$-graded simple current extension
  $$
    \W(c_m)
    := L(c_m,0)\oplus L(c_m,h^{(m)}_{m+1,1})
  $$
  has a unique simple rational vertex operator algebra structure extending $L(c_m,0)$
  if $m\equiv 0,\,3 \pmod{4}$, and a unique simple rational vertex operator
  superalgebra structure extending $L(c_m,0)$ if $m\equiv 1,\,2 \pmod{4}$.
\end{thm}

By this theorem, we introduce the following notion.

\begin{df}\label{df:2.10}
  Let $m\equiv 0$ or $3 \pmod{4}$.
  A simple $c=c_m$ Virasoro vector $u$ of a VOA $V$ is called
  {\it extendable\/} if there exists a non-zero highest weight vector
  $w\in V$ of weight $h^{(m)}_{m+1,1}=m(m+1)/4$ with respect to
  $\vir(u)$ such that the subalgebra generated by $u$ and $w$ is
  isomorphic to the extended Virasoro VOA $\W(c_m)$.
  Note that $w$ is a $\vir(u)$-primary vector, i.e.,
  $L^u(m)w=0 $ for all $m>0$, where $L^u(m)=u_{(m+1)}$.
  We will call such a $w$ an {\it $h_{m+1,1}^{(m)}$-primary vector\/}
  associated to $u$.
\end{df}

\begin{lem}\label{lem:2.11}
  Let $m\equiv 0,\,3 \pmod{4}$ and $u\in V$ be a simple extendable
  $c=c_m$ Virasoro vector.
  Then an $h_{m+1,1}^{(m)}$-primary vector associated to $u$
  is unique up to scalar multiple.
\end{lem}

\pf
Let $w$, $w'$ be $h_{m+1,1}^{(m)}$-primary vectors associated to $u$
and $\W$ the subalgebra generated by $u$ and $w$.
Then the $\W$-submodule $\W'$ generated by $w'$ is isomorphic to
the adjoint module $\W$.
Since we have assumed that $V_0=\C \vac$, we see that $\W =\W'$ and
the assertion follows.
\qed

\medskip

Next we discuss the irreducible modules for $\W(c_m)$ when $m\equiv 0$ or $3 \pmod{4}$.

\begin{lem}[\cite{LLY}]\label{lem:2.11a}
Let $L(c_m,h_{m+1,1}^{(m)})\times L(c_m,h_{r,s}^{(m)}) =L(c_m,\tilde{h}_{r,s}^{(m)})$,
where
\[
\tilde{h}^{(m)}_{r,s}=
\begin{cases}
h^{(m)}_{m+2-r,s} & \text{ when~ } m\equiv 0\pmod{4}, 
\vsb\\
h^{(m)}_{r,m+3-s} & \text{ when~ } m\equiv 3\pmod{4},
\end{cases}
\]
by Eq.\ \eqref{eq:2.2}.
Then $\Delta =h^{(m)}_{r,s}-\tilde{h}^{(m)}_{r,s}$ is in either $\Z$ or $\frac{1}{2}+\Z$.
More precisely, one has:
\\
(1) When $m\equiv 0 \pmod{4}$, $\Delta \in \Z$ if $r$ is odd, 
and $\Delta \in \frac{1}2 +\Z$ if $r$ is even.
\\
(2) When $m\equiv 3 \pmod{4}$, $\Delta \in \Z$ if $s$ is odd,  
and $\Delta \in \frac{1}2 +\Z$ if $s$ is even.
\\
(3) $\Delta =0$, i.e.,  $h^{(m)}_{r,s}= \tilde{h}^{(m)}_{r,s}$ when the triple $(m,r,s)$ satisfies  
$r=(m+2)/2$ if $m\equiv 0\pmod{4}$ or $s=(m+3)/2$ if $m\equiv 3 \pmod{4}$.
\end{lem}

\begin{thm}[\cite{LLY}]\label{thm:2.12a}
Suppose $m\equiv 0$ or $3 \pmod{4}$.
Let $M=L(c_m,h_{r,s}^{(m)})$  be an irreducible $L(c_m,0)$-module and
$\tilde{M}=L(c_m,\tilde{h}^{(m)}_{r,s}):=L(c_m,h^{(m)}_{m+1,1})\fusion M$.
\\
(1) If $h^{(m)}_{r,s} -\tilde{h}^{(m)}_{r,s}\in \Z\setminus \{0\}$, then $M\oplus
\tilde{M}$ affords a unique
structure of an irreducible (untwisted) $\W(c_m)$-module extending $M$.
\\
(2)  If $h^{(m)}_{r,s} -\tilde{h}^{(m)}_{r,s}\in \frac{1}2 +\Z$, then $M\oplus
\tilde{M}$ affords a unique structure of an irreducible $\Z_2$-twisted $\W(c_m)$-module
extending $M$.
\\
(3) If $h^{(m)}_{r,s}= \tilde{h}^{(m)}_{r,s}$,  then $M\oplus \tilde{M}$ is a direct
sum of two inequivalent irreducible (untwisted) $\W(c_m)$-modules. In this
case, there exists two inequivalent structures of an irreducible (untwisted)
$\W(c_m)$-module on $M$ and these structures are $\Z_2$-conjugates of each other. 
We denote them by $M^\pm$.
\end{thm}

Note that Theorem \ref{thm:2.12a} together with Lemma \ref{lem:2.11a} provides a 
classification of irreducible $\W(c_m)$-modules.

\begin{rem}\label{rem:2.10}
  If $u$ is a simple extendable $c=c_m$ Virasoro vector of $V$,
  then it follows from Lemma~\ref{lem:2.11a} and Theorem~\ref{thm:2.12a} that the automorphism
$\tau_u$
  defined in Theorem \ref{thm:2.2} is trivial on $V$ since $V$ is an untwisted
  module over the extended subalgebra $\W(c_m)$ of $\vir(u)$.
\end{rem}

\medskip

In this paper, we will frequently consider simple extendable Virasoro vectors
with central charges $c_3=4/5$ and $c_4=6/7$.
The key feature is that the extended Virasoro VOAs $\W(\sfr{4}{5})$ and
$\W(\sfr{6}7)$ have some natural $\Z_3$-symmetries among their irreducible modules.

\medskip

The extended Virasoro VOA $\W(\sfr{4}{5})$ is rational and has
six inequivalent irreducible modules (cf.\ \cite{KMY}).
They are of the following forms as $L(\sfr{4}{5},0)$-modules:
$$
  L(\sfr{4}{5},0)\oplus L(\sfr{4}{5},3),~~
  L(\sfr{4}{5},\sfr{2}{5})\oplus L(\sfr{4}{5},\sfr{7}{5}),~~
  L(\sfr{4}{5},\sfr{2}{3})^\pm,~~L(\sfr{4}{5},\sfr{1}{15})^\pm ,
$$
where the ambiguity on choosing signs $\pm$ is solved by fusion rules
(cf.\ \cite{M2}).
The extended Virasoro VOA $\W(\sfr{6}{7})$ is rational and has
nine inequivalent irreducible modules (cf.\ \cite{LLY,LY2}).
They are of the following forms as $L(\sfr{6}{7},0)$-modules:
$$
\begin{array}{l}
  L(\sfr{6}{7},0)\oplus L(\sfr{6}{7},5),~~
  L(\sfr{6}{7},\sfr{1}{7})\oplus L(\sfr{6}{7},\sfr{22}{7}),~~
  L(\sfr{6}{7},\sfr{5}{7})\oplus L(\sfr{6}{7},\sfr{12}{7}),
  \vsb\\
  L(\sfr{6}{7},\sfr{4}{3})^\pm,~~~~~~~~~~~
  L(\sfr{6}{7},\sfr{1}{21})^\pm,~~~~~~~~~~~
  L(\sfr{6}{7},\sfr{10}{21})^\pm,
\end{array}
$$
where the ambiguity on choosing signs $\pm$ is again solved by fusion rules
(cf.\ loc.\ cit.).
The fusion rules among irreducible $\W(\sfr{4}{5})$-modules
and $\W(\sfr{6}{7})$-modules are computed in \cite{M2,LLY,LY2} and
they have some natural $\Z_3$-symmetries.
We can extend these symmetries to automorphisms of VOAs containing these extended
Virasoro VOAs as follows.

\begin{thm}[\cite{M2,LLY,LY2}]\label{thm:2.12}
  Let $V$ be a VOA and let $U$ be a sub-VOA of V.
  \\
  (1) Suppose $U\simeq \W(\sfr{4}{5})$.
  Define a linear automorphism $\xi_U$ of $V$ to act on each
  irreducible $\W(\sfr{4}{5})$-submodule $M$ by
  $$
  \begin{cases}
    ~1 & \mbox{if}~~ M\simeq L(\sfr{4}{5},0)\oplus L(\sfr{4}{5},3)
    ~~\mbox{or}~~ L(\sfr{4}{5},\sfr{2}{5})\oplus L(\sfr{4}{5},\sfr{7}{5}),
    \vsb\\
    ~e^{\pm 2\pi\sqrt{-1}/3} & \mbox{if}~~
    M\simeq L(\sfr{4}{5},\sfr{2}{3})^\pm ~~\mbox{or}~~
    L(\sfr{4}{5},\sfr{1}{15})^\pm
  \end{cases}
  $$
  as $L(\sfr{4}{5},0)$-modules.
  Then $\xi_U$ defines an element in $\aut(V)$ satisfying ${\xi_U}^3=1$.
  \\
  (2) Suppose $U\simeq \W(\sfr{6}{7})$.
  Define a linear automorphism $\xi_U$ of $V$ to act
  on each irreducible $\W(\sfr{6}{7})$-submodule $M$ by
  $$
  \begin{cases}
    \, 1 & \mbox{if}~ M\simeq L(\sfr{6}{7},0)\oplus L(\sfr{6}{7},5),~
    L(\sfr{6}{7},\sfr{1}{7})\oplus L(\sfr{6}{7},\sfr{22}{7})~
    \mbox{or}~
    L(\sfr{6}{7},\sfr{5}{7})\oplus L(\sfr{6}{7},\sfr{12}{7}),
    \vsb\\
    \, e^{\pm 2\pi\sqrt{-1}/3} & \mbox{if}~
    M\simeq
    L(\sfr{6}{7},\sfr{4}{3})^\pm,~
    L(\sfr{6}{7},\sfr{1}{21})^\pm~
    \mbox{or}~L(\sfr{6}{7},\sfr{10}{21})^\pm
  \end{cases}
  $$
  as $L(\sfr{6}{7},0)$-modules.
  Then $\xi_U$ defines an element in $\aut(V)$ satisfying ${\xi_U}^3=1$.
\end{thm}

\medskip

\begin{rem}\label{rem:2.13}
If a simple $c=6/7$ Virasoro vector $x\in V$ is extendable,
then $x$ is of $\sigma$-type if and only if the automorphism $\xi_U$
defined in (2) of Theorem \ref{thm:2.12} is trivial on $V$, where
$U$ is the subalgebra isomorphic to $\W(\sfr{6}{7})$ generated by $x$ and
its $5$-primary vector.
\end{rem}


\section{Commutant subalgebras associated to root lattices}\label{sec:4}

In this section, we will construct sub-VOAs of the lattice VOA
$V_{\sqrt{2}{E_6}}$ which will correspond to dihedral subgroups of
the largest Fischer group.
Our construction is similar to the construction in \cite{LYY1} and \cite{HLY} in
the case of the root lattices $E_8$ and $E_7$.

\smallskip




\subsection{The algebras $\UFnX{nX}$ and $\VFnX{nX}$}\label{sec:4.1}

\paragraph{Commutant subalgebras.}

Let $(V,\w)$ be a VOA and $(U,e)$ be a sub-VOA.
Then the commutant subalgebra of $U$ is defined by
\begin{equation}\label{eq:3.1}
  \com_V(U):= \{ a\in V \mid a_{(n)} U=0 \text{ for all } n\geq 0\} .
\end{equation}
It is known (cf.\ \cite{FZ}) that
\begin{equation}\label{eq:3.2}
  \com_V(U)=\ker_V e_{(0)}
\end{equation}
and in particular $\com_V(U)=\com_V(\vir(e))$. Therefore, the commutant
subalgebra of $U$ is determined only by the conformal vector $e$ of $U$. It is
also shown in Theorem 5.1 of \cite{FZ} that $\w -e$ is also a Virasoro vector
if $\w_{(2)} e=0$. Provided $V_1=0$, we always have $\w_{(2)}e=0$. In that case,
we have two mutually commuting subalgebras $\com_V(\vir(e))=\ker_V e_{(0)}$ and
$\com_V(\vir(\w-e))=\ker_V(\w-e)_{(0)}$ and the tensor product $\com_V
(\vir(\w-e))\tensor \com_V(\vir(e))$ forms an extension of $\vir(e)\tensor
\vir(\w-e)$. More generally, we say a sum $\w =e^1+\cds +e^n$ is a {\it
Virasoro frame\/} if all $e^i$ are Virasoro vectors and
$[Y(e^i,z_1),Y(e^j,z_2)]=0$ for $i\ne j$.

\paragraph{The algebras $\UFnX{nX}$.}

Let $\alpha_1$,~$\ldots$,~$\alpha_6$ be a system of simple roots for $E_6$.
We let $\alpha_0$ be the root such that
$-\alpha_0=\sum_{i=1}^6 m_i \alpha_i$ is the highest root for the chosen
simple roots.
Note that all $m_i$ are positive integers.
We also set $m_0=1$.
For any $i=0$, $\ldots$, $6$, we consider the sublattice $L_i$ of $E_6$
generated by the roots $\alpha_j$, $0\leq j\leq 6$, $j\not=i$.
One observes that $L_i$ is also of rank $6$ and the quotient group
$E_6/L_i$ is cyclic of order $m_i$ with generator $\alpha_i+L_i$.
Thus one has
\begin{equation}\label{eq:4.1}
  E_6 = L_i \sqcup (\al_i+L_i) \sqcup (2\alpha_i+L_i)\sqcup
  \cdots \sqcup ((m_i-1)\al_i+L_i).
\end{equation}

We denote by $R_1$, $\dots$, $R_\ell$ the indecomposable components of the
lattice $L_i$ which are root lattices of type $A_n$, $D_n$ or $E_6$. Hence
$L_i=R_1\oplus \cdots \oplus R_\ell$ where the direct sum of lattices denotes
the orthogonal sum. In fact, the Dynkin diagram of $L_i$ is obtained from the
affine Dynkin diagram of $E_6$ by removing the node~$\alpha_i$ and the adjacent
edges. We recall here that the affine Dynkin diagram of $E_6$ is the graph with
vertex set $\{\alpha_0, \ldots,  \alpha_6\}$ and two nodes $\alpha_i$ and
$\alpha_j$, $0\leq i$, $j\leq 6$, are joined by an edge if $\la \alpha_i,
\alpha_j \ra=-1$. The diagram has the following form:
\begin{equation}\label{eq:4.21}
\begin{array}{l}
  \hspace{120.5pt}\circ \hspace{12pt}  \alpha_0\vspace{-14pt}\\
  \hspace{121.3pt}| \vspace{-10pt}\\
  \hspace{121.3pt}| \vspace{-10pt}\\
  \hspace{121.3pt}| \vspace{-14pt}\\
  \hspace{120.5pt}\circ \hspace{12pt}  \alpha_6
  \vspace{-14pt}\\
  \hspace{121.3pt}| \vspace{-10pt}\\
  \hspace{121.3pt}| \vspace{-10pt}\\
  \hspace{121.3pt}| \vspace{-14pt}\\
  \hspace{8pt}
  \circ\hspace{-5pt}-\hspace{-7pt}-\hspace{-7pt}-\hspace{-7pt}-
\hspace{-7pt}-\hspace{-7pt}\hspace{-7pt}-\hspace{-7pt}-\hspace{-7pt}-\hspace{-7pt}-
\hspace{-7pt}-\hspace{-6pt}-\hspace{-5pt}
\circ\hspace{-5pt}-\hspace{-7pt}-\hspace{-7pt}-\hspace{-7pt}-
\hspace{-7pt}-\hspace{-7pt}\hspace{-7pt}-\hspace{-7pt}-\hspace{-6pt}-
\hspace{-6pt}-\hspace{-6pt}-\hspace{-6pt}-\hspace{-5pt}
\circ\hspace{-5pt}-\hspace{-7pt}-\hspace{-7pt}-\hspace{-7pt}-
\hspace{-7pt}-\hspace{-7pt}\hspace{-7pt}-\hspace{-7pt}-\hspace{-6pt}-
\hspace{-6pt}-\hspace{-6pt}-\hspace{-6pt}-\hspace{-5pt}
\circ\hspace{-5pt}-\hspace{-7pt}-\hspace{-7pt}-\hspace{-7pt}-
\hspace{-7pt}-\hspace{-7pt}\hspace{-7pt}-\hspace{-7pt}-\hspace{-6pt}-
\hspace{-6pt}-\hspace{-6pt}-\hspace{-6pt}-\hspace{-5pt}\circ
  \vspace{-5pt}\\
  \hspace{7pt} \alpha_1
  \hspace{42pt} \alpha_2
  \hspace{45pt} \alpha_3
  \hspace{45pt} \alpha_4
  \hspace{45pt} \alpha_5
  \end{array}
\end{equation}

\smallskip

The decomposition \eqref{eq:4.1} of the lattice $E_6$ leads to the
decomposition
$$
  V_{\sqrt{2}E_6}=\bigoplus_{r=0}^{m_i-1}V_{\sqrt{2}(r\alpha_i+L_i)}
$$
of the lattice VOA $V_{\sqrt{2}E_6}$.
We define a linear map $\rho_i: V_{\sqrt{2}E_6} \rightarrow V_{\sqrt{2}E_6}$
by
\begin{equation}\label{eq:4.2}
  \rho_{i}(u)
  = \zeta_{m_i}^r u \quad \text{ for } u \in  V_{\sqrt{2}(r \al_i+L_i)},
  ~~~~~\mbox{where}~~ \zeta_{m_i} =e^{2 \pi \sqrt{-1}/m_i}.
\end{equation}
Then $\rho_i$ is an element of $\aut (V_{\sqrt{2}E_6})$ of order $m_i$
and the fixed point sub-VOA $V_{\sqrt{2}E_6}^{\la\rho_i\ra}$ is exactly
$V_{\sqrt{2}L_i}$.

For a root lattice $S$, we denote by $\Phi(S)$ its root system.
Then, by \cite{DLMN}, the conformal vector $\omega_S$ of $V_{\sqrt{2}S}$
is given by
$$
  \omega_S = \frac{1}{4h}\sum_{\al\in\Phi(S)}\al(-1)^2 \vac,
$$
where $h$ is the Coxeter number of $S$.
Now define
\begin{equation}\label{eq:4.3}
  \tilde{\w}_S
  := \dfrac{2}{h+2}\,\omega_S+\dfrac{1}{h+2}
     \dsum_{\al\in\Phi(S)}e^{\sqrt{2}\alpha}.
\end{equation}
It is shown in~\cite{DLMN}
that $\tilde{\w}_S$ is a Virasoro vector of central charge
$2n/(n+3)$ if $S$ is of type $A_n$, $1$ if $S$ is of type $D_n$, and
$6/7$, $7/10$ and $1/2$ if $S$ is of type $E_6$, $E_7$, $E_8$, respectively.
From the irreducible decomposition
$L_i=R_1\oplus \cds \oplus R_\ell \subset E_6$,
we have sublattices $R_s$ of $E_6$ and obtain a factorization
\begin{equation}\label{eq:4.4}
  V_{\sqrt{2}L_i}= V_{\sqrt{2}R_1}\tensor \cds \tensor V_{\sqrt{2}R_\ell}
  \subset V_{\sqrt{2}E_6}.
\end{equation}
Associated to the root subsystems $\Phi(R_s)$ of $\Phi(E_6)$, we also have
Virasoro vectors
\begin{equation}\label{eq:4.5}
  \w^s
  :=\tilde{\w}_{R_s}
  = \dfrac{2}{h_s+2}\,\omega_{R_s}+\dfrac{1}{h_s+2}
    \dsum_{\al\in\Phi(R_s)}e^{\sqrt{2}\alpha}\in V_{\sqrt{2}R_s}
    \subset V_{\sqrt{2}E_6},
    \quad 1\leq s\leq \ell,
\end{equation}
where $\w_{R_s}$ is the conformal vector of $V_{\sqrt{2}R_s}$ and $h_s$ is the
Coxeter number of $R_s$. It follows from the definition that $\w^s, 1\leq s\leq
\ell,$ are mutually orthogonal simple Virasoro vectors in $V_{\sqrt{2}E_6}$.
Consider
$$
  X^r:=\sum_{\scriptstyle \beta \in r\alpha_i+ L_i  \atop \scriptstyle
  \la \beta, \beta \ra=2} e^{\sqrt{2}\beta},
  \qquad 1\leq r\leq m_i-1,
$$
in the weight two subspace of $V_{\sqrt{2}E_6}$.
It is shown in Proposition 2.2 of~\cite{LYY1} that the vectors $X^r$ are
highest weight vectors for $\vir(\w^1)\tensor \cds \tensor \vir(\w^\ell)$
with total weight~$2$.

\smallskip

Since all $\w^s$, $1\leq s\leq \ell$, are contained in the fixed point sub-VOA
$V_{\sqrt{2}E_6}^+$, which has a trivial weight one subspace,
the vector $\w_{E_6}-(\w^1+\cds +\w^\ell)$ is a Virasoro vector of
$V_{\sqrt{2}E_6}$ as discussed at the beginning of the section.
We are interested in the following commutant subalgebras:
\begin{df}\label{UnX}
For $i=0, \dots, 6$,   let $L_i < E_6$ be defined as in \eqref{eq:4.1} and $R_1$,
$\dots$, $R_\ell$ the indecomposable components of $L_i$.  Let
$\w^s=\tilde{\w}_{R_s}$ be Virasoro vectors defined as in \eqref{eq:4.5} for
$1\leq s\leq \ell$. The algebra $U(i)$ is the vertex operator algebra
$$\begin{array}{lrl}
  U(i)
  &=&\mathrm{Com}_{V_{\sqrt{2}E_6}}( \vir(\w_{E_6}-(\w^1+\cds +\w^\ell)))
  \vsb\\
  &=&\ker_{V_{\sqrt{2}E_6}} \l(\w_{E_6}-(\w^1+\cds +\w^\ell)\r)_{(0)}.
\end{array}$$
 \end{df}
It is clear from the construction that $U(i)$ has a Virasoro frame
$\w^1+\cds +\w^\ell$.
We will consider an embedding of $U(i)$ into a larger VOA and then describe the
commutant algebra $U(i)$ using the larger VOA.

\medskip

It is clear that $U(i)$ forms an extension of 
$\vir(\w^1)\tensor \cds \tensor \vir(\w^\ell)$ and contains highest weight
vectors $X^r$, $1\leq r<  m_i$. We will see in Section \ref{sec:5} that we can
embed $U(i)$ into the Moonshine VOA and therefore $U(i)$ has a trivial weight
one subspace. Consequently, the weight two subspace of $U(i)$ carries a
structure of a commutative non-associative algebra called the Griess algebra of
$U(i)$, even though $V_{\sqrt{2}E_6}$ has a non-trivial weight one subspace.
In Section \ref{sec:3.3new}, we will explicitly describe the Griess algebra of $U(i)$. Namely, we will show that the Griess algebra of $U(i)$ is given by
$$
  \mathcal{G}(i):
  = \Span_\C \{~ \w^s,~ X^r \mid 1\leq s\leq \ell,~~
     1\leq r\leq m_i-1\},
$$
which is of dimension $\ell+m_i-1$.

Recall $\rho_i\in \aut(V_{\sqrt{2}E_6})$ defined as in \eqref{eq:4.2}.
By definition, it is clear that $\tilde{\w}_{E_6}$ and $\rho_i \tilde{\w}_{E_6}$ are
linear combinations of $\om^s$ and $X^r$, and hence are contained in
$\mathcal{G}(i)\subset U(i)$. We will also discuss the structure of the subalgebra
generated by $\tilde{\w}_{E_6}$ and $\rho_i\tilde{\w}_{E_6}$ and compare it with
$\mathcal{G}(i)$ in Section~\ref{sec:3.3new}.

\medskip


\paragraph{The algebras $\VFnX{nX}$.}

We also consider another class of commutant algebras inside the VOA
$V_{\sqrt{2}E_8}$. These commutant algebras will be used in Section~\ref{sec:5.2.3}
to show that $U(i)$ defined in Definition~\ref{UnX} can be embedded into the Fischer group VOA $\VF$.

\smallskip

We fix an embedding of $E_6$ into $E_8$. Let
\begin{equation}\label{eq:4.7}
  Q:= \mathrm{Ann}_{E_8}(E_6)
  = \{\al\in E_8\mid \la \al, E_6\ra =0\}.
\end{equation}
Then $Q\simeq A_2$ and $Q\oplus E_6$ forms a full rank sublattice of $E_8$.
Note that such an embedding is unique up to an automorphism of $E_8$.

Recall that $L_i$ is the sublattice of $E_6$ generated by roots
$\alpha_j$, $j\ne i$.
Then we have an embedding of $\EL_i:=Q \oplus L_i$ into $E_8$.
Since $L_i$ is a full rank sublattice of $E_6$, $\EL_i$ is also a full rank
sublattice of $E_8$.
Thus $E_8/\EL_i$ is a finite abelian group whose order is $3m_i$.
We fix the corresponding embedding
$V_{\sqrt{2}\EL_i} \subset V_{\sqrt{2}E_8}$.

We have the decomposition
$\EL_i =Q \oplus R_1\oplus \cdots \oplus R_{\ell}$
into a sum of irreducible root lattices which gives rise to a factorization
$$
  V_{\sqrt{2} \EL_i}
  = V_{\sqrt{2} Q} \tensor V_{\sqrt{2}L_i}
  = V_{\sqrt{2} Q}\tensor V_{\sqrt{2} R_1}\tensor \cds \tensor
  V_{\sqrt{2} R_\ell} \subset V_{\sqrt{2}E_8}.
$$
Let $\w_{E_8}$ be the conformal vector of $V_{\sqrt{2}E_8}$ and
let $\tilde{\om}_{Q}\in V_{\sqrt{2}Q}$ and $\om^s\in V_{\sqrt{2}R_s}$ be
the Virasoro vectors defined as in \eqref{eq:4.3} and \eqref{eq:4.5}, respectively.
By the same argument as for $U(i)$, one sees
$\w_{E_8}-(\tilde{\w}_{Q} +\w^1+ \cdots+\w^\ell)$ is  a Virasoro vector of
$V_{\sqrt{2}E_8}$, and we can define a commutant subalgebra:
\begin{df}\label{VnX}
The algebra $V(i)$ is the commutant subalgebra
$$  V(i):=\com_{V_{\sqrt{2}E_8}}( \vir(\om_{E_8}-(\tilde{\om}_{Q} +{\om}^1
  +\cdots +\om^\ell))). $$
\end{df}

\begin{rem}\label{rem:4.1}
  As explained in~\cite{HLY}, $\com_{V(i)}(\vir(\tilde{\om}_Q))$ coincides
  with $U(i)$.
\end{rem}

We finally set
\begin{equation}\label{eq:4.9}
  \mathcal{F}(i)
  := \{ g\in \aut(V_{\sqrt{2}E_8}) \mid g=\id \ \text{ on }\
     V_{\sqrt{2}\EL_{i}} \} .
\end{equation}
Then $\mathcal{F}(i)$ is canonically isomorphic to the group of characters
of $E_8/\EL_{i}$.
The subalgebra $V(i)$ of $V_{\sqrt{2}E_8}$ is invariant under the action of
$\mathcal{F}(i)$  since all $\tilde{\om}_Q$, $\om^1$, $\dots$,
$\om^\ell$ and the conformal vector $\w_{E_8}$ of $V_{\sqrt{2}E_8}$
are clearly fixed by $\mathcal{F}(i)$.
Note that the special Ising vector
\begin{equation}\label{eq:4.10}
  \hat e
  :=\tilde{\om}_{E_8}
  = \dfr{1}{16}\w_{E_8} +\dfr{1}{32}
  \sum_{\alpha\in \Phi(E_8)} e^{\sqrt{2}\alpha} \in V_{\sqrt{2}E_8}
\end{equation}
is contained in $V(i)$ (cf.\ \cite{LYY1,LYY2}) and thus
$\{ g \hat{e} \in V_{\sqrt{2}E_8} \mid g\in \mathcal{F}(i)\} \subset V(i)$.

\medskip

%
%
%

\paragraph{McKay's $E_6$-correspondence.}

We like to explain McKay's correspondence between the conjugacy
classes 1A, 2A and 3A of the Fischer group $\mathrm{Fi}_{24}$
which are the products of 2C-involutions of the Fischer group and
the numerical labels $m_i$ of the affine $E_6$ Dynkin diagram
as given by the following figure:
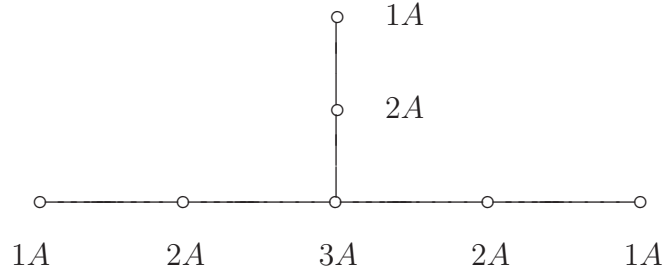
\begin{equation*}\label{eq:4.22}
\begin{array}{l}
  \hspace{120.5pt}\circ \hspace{15pt}  1A\vspace{-14pt}\\
  \hspace{121.3pt}| \vspace{-14pt}\\
  \hspace{121.3pt}| \vspace{-10pt}\\
  \hspace{121.3pt}| \vspace{-14pt}\\
  \hspace{120.5pt}\circ \hspace{15pt}  2A
  \vspace{-14pt}\\
  \hspace{121.3pt}| \vspace{-14pt}\\
  \hspace{121.3pt}| \vspace{-10pt}\\
  \hspace{121.3pt}| \vspace{-14pt}\\
  \hspace{8pt}
  \circ\hspace{-5pt}-\hspace{-7pt}-\hspace{-7pt}-\hspace{-7pt}-
  \hspace{-7pt}-\hspace{-7pt}\hspace{-7pt}-\hspace{-7pt}-\hspace{-7pt}-
  \hspace{-7pt}-\hspace{-7pt}-\hspace{-6pt}-\hspace{-5pt}\circ
  \hspace{-5pt}-\hspace{-7pt}-\hspace{-7pt}-\hspace{-7pt}-
  \hspace{-7pt}-\hspace{-7pt}\hspace{-7pt}-\hspace{-7pt}-\hspace{-6pt}-
  \hspace{-6pt}-\hspace{-6pt}-\hspace{-6pt}-\hspace{-5pt}
  \circ\hspace{-5pt}-\hspace{-7pt}-\hspace{-7pt}-\hspace{-7pt}-
  \hspace{-7pt}-\hspace{-7pt}\hspace{-7pt}-\hspace{-7pt}-\hspace{-6pt}-
  \hspace{-6pt}-\hspace{-6pt}-\hspace{-6pt}-\hspace{-5pt}
  \circ\hspace{-5pt}-\hspace{-7pt}-\hspace{-7pt}-\hspace{-7pt}-
  \hspace{-7pt}-\hspace{-7pt}\hspace{-7pt}-\hspace{-7pt}-\hspace{-6pt}-
  \hspace{-6pt}-\hspace{-6pt}-\hspace{-6pt}-\hspace{-5pt} \circ
  \vspace{-1pt}\\
  1A\hspace{44pt}  2A \hspace{43  pt} 3A\hspace{43pt} 2A\hspace{43pt} 1A
  \end{array}
\end{equation*}
Note that the correspondence is not one-to-one but only
up to the diagram automorphism.

\medskip

Because of this correspondence, we change our notation slightly
and denote $L_i$ by $L_{nX}$, $\rho_i$ by $\rho_{nX}$, $\EL_i$ by $\EL_{nX}$,
$\mathcal{F}(i)$ by $\mathcal{F}_{nX}$, $V(i)$ by $\VFnX{nX}$, $U(i)$ by $\UFnX{nX}$
and $\mathcal{G}(i)$ by $\GFnX{nX}$, where
$nX\in\{1A,\, 2A,\, 3A \}$ is the label of the corresponding node
in \eqref{eq:4.22}.
Explicitly, we have:
\begin{equation}\label{eq:4.23}
  L_{1A}\simeq E_6,\qquad
  L_{2A}\simeq A_1\oplus A_5,\qquad
  L_{3A}\simeq A_2\oplus A_2 \oplus A_2.
\end{equation}
We also have that  
$$
  \tilde{\om}_Q
  = \frac{1}{15}\left( {\beta_{0}}(-1)^{2}+{\beta_{1}}(-1)^{2}
  +{\beta_{2}}(-1)^2\right) \vac
  +\frac{1}{5}\sum_{\al\in \Phi(Q)} e^{\sqrt{2}\al}
  \in V_{\sqrt{2}Q}
$$
is a simple $c=4/5$ Virasoro vector in $V_{\sqrt{2}Q}$,
where $\{\beta_1, \beta_2\}$ is a set of simple roots for
$Q \simeq A_2$ and $\beta_0=-(\beta_1+\beta_2)$.


\subsection{Structures of  $\VFnX{nX}$ and $\UFnX{nX}$.}\label{sec:3.3new}
We determine the structures of $\VFnX{nX}$ and $\UFnX{nX}$.

\paragraph{\bf 1A case.}
In this case, we have  $\EL_{1A}\simeq A_2\oplus E_6$ and
 $\VFnX{1A}\simeq U_{3A}$ by \cite{LYY2}.
It is clear that the following holds
(see \cite{LY2} and Section~\ref{sec:3} for details):
\begin{lem}\label{U1A}
$\UFnX{1A}\simeq \W(\sfr{6}{7})=L(\sfr{6}{7},0)\oplus L(\sfr{6}{7},5)$.
\end{lem}
Therefore the weight two subspace of $\UFnX{1A}$ is one-dimensional.

\paragraph{\bf 2A case.}
In this case $\EL_{2A}\simeq A_2\oplus A_1\oplus A_5$ and
$E_8/\EL_{3A}\simeq \Z_6$.
By construction, the commutant subalgebra $\VFnX{2A}$ is the same as
the monstrous 6A-algebra $U_{6A}$ discussed in \cite{LYY2}
(see also the Appendix of \cite{HLY}).
Thus, the following result follows from \cite{LYY2}:
\begin{lem}\label{U2A}
There is the decomposition
$$\begin{array}{ll}
\UFnX{2A}\ \simeq\
   &\W(\sfr{6}{7})\otimes  L(\sfr{25}{28},0) \ \oplus \
  \bigl(L(\sfr{6}{7}, \sfr{5}{7})\oplus L(\sfr{6}{7},\sfr{12}{7})\bigr)
    \otimes L(\sfr{25}{28},\sfr{9}7)\qquad\qquad
\\
   &\ \oplus \  \bigl(L(\sfr{6}{7},\sfr{1}{7})
    \oplus L(\sfr{6}{7},\sfr{22}{7})\bigr)
    \otimes L(\sfr{25}{28},\sfr{34}{7})
\end{array}
$$
as a $\W(\sfr{6}{7})\otimes  L(\sfr{25}{28},0)$-module.
\end{lem}
It follows from the decomposition in the Lemma that the weight two subspace of
$\UFnX{2A}$ is $3$-dimensional and coincides with $\GFnX{2A}$.

\paragraph{\bf 3A case.}
In this case, $\EL_{3A}\simeq A_2\oplus A_2\oplus
A_2\oplus A_2$ and $E_8/\EL_{3A}\simeq\Z_3\oplus \Z_3$.
In fact, the coset structure $E_8/\EL_{3A}$ can be identified with the
ternary {\it tetra code} $\mathcal{C}_4$, whose generator matrix is given by
$$
\begin{bmatrix}
  1 & 1 & 1 & 0
  \\
  1 & -1 & 0 & 1
\end{bmatrix}_.
$$
The sub-VOA $\VFnX{3A}$ is indeed the ternary code VOA defined in \cite{KMY}:
\begin{equation}\label{eq:3.13+}
  M_{\mathcal{C}_4}
  = \bigoplus_{\alpha=(\alpha_1,\alpha_2,\alpha_3,\alpha_4)\in \mathcal{C}_4}
  L_{4/5}(\alpha_1)\tensor L_{4/5}(\alpha_2)\tensor L_{4/5}(\alpha_3)\tensor L_{4/5}(\alpha_4),
\end{equation}
where $L_{4/5}(0)=L(\sfr{4}{5},0)\oplus L(\sfr{4}{5},3)=\W(\sfr{4}{5})$ and
$L_{4/5}(\pm 1)=L(\sfr{4}{5},\sfr{2}{3})^\pm$ as $\W(\sfr{4}{5})$-modules.
Note also that the cosets of $L_{3A}$ in $E_6$ can be parameterized by
the ternary repetition code of length~$3$:
\begin{equation}\label{eq:3.14+}
  \mathcal{D}= \{(0,0,0),\, (1,1,1),\, (-1,-1,-1)\}.
\end{equation}
Thus the commutant subalgebra $\UFnX{3A}$ is the ternary code VOA
$M_\mathcal{D}$ defined in \cite{KMY}:
\begin{equation}\label{eq:3.15+}
  M_\mathcal{D}=\l( L(\sfr{4}{5},0)\oplus L(\sfr{4}{5},3)\r)^{\tensor 3}
  \oplus \l( L(\sfr{4}{5},\sfr{2}{3})^+\r)^{\tensor 3}
  \oplus \l( L(\sfr{4}{5},\sfr{2}{3})^-\r)^{\tensor 3}.
\end{equation}
It is shown in \cite{KMY} that
$$
  \ker_{V_{\sqrt{2}A_2}}( \om_{A_2} - \tilde{\om}_{A_2})_{(0)}
  \simeq  \W(\sfr{4}{5})=L(\sfr{4}{5},0)\oplus L(\sfr{4}{5},3)
$$
and we get:
\begin{lem}\label{U3A}
One has a decomposition
$$
  \UFnX{3A} \  \simeq \
 \W(\sfr{4}{5})^{\otimes 3}
 \ \oplus \
 \bigl(L(\sfr{4}{5},\sfr{2}{3})^+\bigr)^{ \otimes 3}
  \ \oplus \
 \bigl( L(\sfr{4}5,\sfr{2}3)^-\bigr)^{ \otimes 3}
$$
as a $\W(\sfr{4}{5})^{\otimes 3}$-module.
\end{lem}
Now we see that the weight two subspace of $\UFnX{3A}$ is $5$-dimensional
and coincides with $\GFnX{3A}$.

\begin{rem} \label{rem:3.7}
By the comments after Remark \ref{rem:4.1}, we know that $\mathcal{F}_{3A}$ acts on
$\VFnX{3A}$. In this case, the fixed point space is
$\VFnX{3A}^{\mathcal{F}_{3A}}\simeq L_{4/5}(0)\tensor L_{4/5}(0)\tensor
L_{4/5}(0)\tensor L_{4/5}(0)\simeq \W(\sfr{4}{5})^{\tensor 4}$  while
$L_{4/5}(\alpha_1)\tensor L_{4/5}(\alpha_2)\tensor L_{4/5}(\alpha_3)\tensor L_{4/5}(\alpha_4)$,
$(\alpha_1,\alpha_2,\alpha_3,\alpha_4)\in \mathcal{C}_4$, are character spaces of
$\mathcal{F}_{3A}$ in $\VFnX{3A}$.  The automorphism group of $\VFnX{3A}\simeq M_{\mathcal{C}_4}$
was computed in \cite[Proposition 5.4]{KMY}. It is isomorphic to $3^2.(2\mathrm{S}_4)\simeq
3^2.\aut(\mathcal{C}_4)\simeq \mathrm{AGL}_2(3)$.
The subgroup $\mathcal{F}_{3A}: \la \tau_{\hat{e}}\ra \simeq 3^2:2$ fixes
the four $c=4/5$ Virasoro vectors and is the stabilizer of the sub-VOA $\VFnX{3A}^{\mathcal{F}_{3A}}\simeq
L_{4/5}(0)\tensor L_{4/5}(0)\tensor L_{4/5}(0)\tensor L_{4/5}(0)\simeq
\W(\sfr{4}{5})^{\tensor 4}$.
\end{rem}

In the Appendix,  we will need a result about generators for~$\VFnX{3A}$.
\begin{lem}\label{VF3A}
Let $\rho_1$ and $\rho_2$ be generating ${\cal F}_{3A}\simeq 3^2$. Then
$\VFnX{3A}$ is generated by $\hat e$, $\rho_1\hat e$ and $\rho_2\hat e$.
\end{lem}

\pf
By Remark~\ref{rem:3.7}, we know that  $\mathcal{F}_{3A}$ acts on $\VFnX{3A} < V_{\sqrt{2}E_8}$ and
the fixed point subspace is $\VFnX{3A}^{\mathcal{F}_{3A}}\simeq
\W(\sfr{4}{5})^{\tensor 4}$.
It is clear that
$e^{i,j}: = \rho_1^i \rho_2^j\hat e$,
$0 \leq i$, $j \leq 2$, are contained in $\VFnX{3A}$ since $\hat e\in \VFnX{3A}$.


Let $W$ be the sub-VOA of $\VFnX{3A}$ generated by $ \{e^{i,j}\mid 0 \leq i,\,j \leq 2\}$.
For each $i$, $j=0$, $1$, $2$ with $(i,j)\neq (0,0)$, the Ising vectors $e^{0,0}$ and $e^{i,j}$
generate a sub-VOA in $\VFnX{3A}$ isomorphic to $U_{3A}$ which contains a subalgebra isomorphic to
$\W(\sfr{4}{5})$ fixed by $\mathcal{F}_{3A}$ (cf.~\cite{LYY2}).
In fact, the $W_3$-algebra $\W(\sfr{4}{5})$ is contained in a lattice sub-VOA
$V_{\sqrt{2}A_2}<V_{\sqrt{2}E_8}$, where $\sqrt{2}A_2<\sqrt{2}A_2^4<\sqrt{2}E_8$.
By varying $i$, $j$ (say, $(i,j)=(1,0)$, $(0,1)$, $(1,1)$ and $(1,2)$),
one can obtain four mutually orthogonal $\W(\sfr{4}{5})$ in
$\VFnX{3A}^{ \mathcal{F}_{3A} }$.
Thus we have $\VFnX{3A}^{ \mathcal{F}_{3A}}\simeq  \W(\sfr{4}{5})^{\tensor 4} < W$.

Now let $\zeta=e^{2\pi i/3}$ and let
\[
  v_{m,n}
  = \sum_{i=0}^2\sum_{j=0}^2 \zeta^{mi+nj} \, e^{-i,-j},
\]
for any $m$, $n=0$, $1$, $2$.  Then by direct calculation, it is easy to verify that
$$
  \rho_1^k \rho_2^\ell (v_{m,n}) = \zeta^{km+\ell n } v_{m,n}.
$$
In other words, $v_{m,n}$ spans a $1$-dimensional
$ \mathcal{F}_{3A}$-submodule affording the character $\chi_{m,n}$,
where $\chi_{m,n}( \rho_1^i  \rho_2^j)= \zeta^{mi+nj}$. It also generates the irreducible
$\VFnX{3A}^{ \mathcal{F}_{3A}}\simeq  \W(\sfr{4}{5})^{\tensor 4}$-submodule
affording the character $\chi_{m,n}$.
Hence, $W$ contains a sub-VOA isomorphic to $M_{\mathcal{C}}\simeq \VFnX{3A}$ and we have
$W= \VFnX{3A}$.

Now note that the 3A-algebra generated by $e^{0,0}$ and $e^{1,0}$ contains
$e^{2,0}= \tau_{e^{0,0}} e^{1,0}$. Similarly, we get $e^{0,2}$.  Since $\rho_1$ and $\rho_2$ are
inverted by $\tau_{\hat{e}}$, we have
\[
\begin{split}
 \tau_{e^{i,0}} ( e^{0,j})  & = \tau_{\rho_1^i \hat{e}} ( \rho_2^j \hat{e})
=
 \rho_1^i \tau_{\hat{e}} \rho_1^{-i} \rho_2^{j} \hat{e}
= \rho_1^{2i} \rho_2^{2j} \hat{e} = e^{2i,2j}
\end{split}
\]
and thus  $e^{2i,2j}$ is contained in the 3A-algebra generated by
$e^{i,0}$ and $e^{0,j}$.  Hence,  $\VFnX{3A}$ is generated by $e^{0,0}=\hat{e}$,
$\rho_1 \hat{e}=e^{1,0}$ and $\rho_2 \hat{e}= e^{0,1}$.
\qed

\begin{rem}\label{32:2}
The irreducible modules for ternary code VOA $M_D$ have been studied in \cite{Lam3}.
It is known (cf.\ Theorem~4.8 and~4.10 of \cite{Lam3}) that if $D$ is a self-dual ternary code,
then all irreducible $M_D$-modules can be realized (using coset or GKO construction)
as submodules of a certain lattice VOA
$V_{\Gamma_D}$. We refer to Section 2 of \cite{KMY} or Section~3.2 of \cite{Lam3}
for the precise definition of $\Gamma_D$. When $D=\mathcal{C}_4$ is the tetra code,
we have $\Gamma_D=\Gamma_{\mathcal{C}_4}\simeq \sqrt{2}E_8$ and all irreducible $M_{\mathcal{C}_4}$-modules
are contained in the lattice VOA $V_{\sqrt{2}E_8}$. The subgroup
$3^2{:}2 \simeq \mathcal{F}_{3A}{:}\la \tau_{\hat{e}}\ra < \aut(V_{\sqrt{2}E_8})$
actually acts faithfully on all irreducible $M_{\mathcal{C}_4}$-submodules.
Thus, if $V$ is a VOA containing $\VFnX{3A}\simeq M_{\mathcal{C}_4}$, then
the involutions $\tau_{\hat e}$, $\tau_{\rho_1\hat e}$ and $\tau_{\rho_2\hat e}$ generate
a group of the shape $3^2{:}2$ in $\aut(V)$.
\end{rem}


\paragraph{Subalgebras generated by $\tv$ and $\tv'$.}
Set
\begin{equation}\label{tv}
  \tv :=\tilde{\w}_{E_6} ~~~~~~~~~\mbox{and}~~~~~~~~~~
  \tv' :=\rho_{nX} \tilde{\w}_{E_6}.
\end{equation}
By definition, the Virasoro vectors $\tv$ and $\tv'$ are contained in
$\UFnX{nX}$. We will discuss $\GFnX{nX}$ and the subalgebras generated by
$\tv$ and $\tv'$.


\paragraph{1A case.}
In this case $L_{1A}\simeq E_6$ and $\rho_{1A}$ is trivial.
Thus $\tv'=\tv$,  $\la \tv,\tv'\ra = 3/7$ and $\tv$ generates $\GFnX{1A}$,
but not $\UFnX{1A}$.

\paragraph{2A case.}
In this case $L_{2A}\simeq A_1\oplus A_5$ and $\ell =2$.

The vectors $\w^1$ and $\w^2$ are Virasoro vectors with central charges
$1/2$ and $5/4$, respectively, and $X=X^1$ is a highest weight vector
for $\vir (\w^1)\tensor \vir (\w^2)$ with highest weight
$(\sfr{1}2,\sfr{3}{2})$. One easily obtains:
\begin{lem}\label{GF2A}
The Griess algebra $\GFnX{2A}$ is spanned by $\om^1$, $\om^2$, and $X$ and we
have the following commutative algebra structure on $\GFnX{2A}$:
$$\begin{array}{c|ccc}
{a}_{(1)}b & \w^1 & \w^2  & X \\ \hline
\w^1       & 2   \w^1  &  0 &  \frac{1}{2}X  \\
\w^2       &    & 2   \w^2 &   \frac{3}{2}X  \\
 X         &    &          &  80\w^1 + 96\w^2
\end{array}
\qquad\qquad\qquad
\begin{array}{c|ccc}
 \la a,b \ra & \w^1 & \w^2  & X \\ \hline
\w^1       & \frac{1}{4}    &  0 &  0  \\
\w^2       &    & \frac{5}{8}  &   0  \\
 X         &    &          & 40
\end{array}
$$
\end{lem}
\medskip

One verifies that
$$
  \tv=\dfr{2}{7}\w^1+\dfr{4}{7}\w^2+\dfr{1}{14} X,~~~~
  \tv'=\dfr{2}{7}\w^1+\dfr{4}{7}\w^2-\dfr{1}{14} X,~~~~
  \la \tv,\tv'\ra = \dfr{1}{49}.
$$
It is also easily verified that $\GFnX{2A}$ is generated by
$\tv$ and $\tv'$.
Set
$$
  u =\dfr{5}{7}\w^1+\dfr{3}{7}\w^2-\dfr{1}{14} X.
$$
Then $\tv$ and $u$ are the mutually orthogonal Virasoro vectors
with central charges $6/7$ and $25/28$, respectively, used in the decomposition
of $\UFnX{2A}$ given before.

\paragraph{3A case.}
In this case $L_{3A}\simeq A_2\oplus A_2\oplus A_2$ and $\ell=3$.

The three vectors $\w^1$, $\w^2$ and $\w^3$ are mutually orthogonal
Virasoro vectors with central charge $4/5$, and $X^1$, $X^2$ are
highest weight vectors  for
$\vir (\w^1)\tensor \vir (\w^2) \tensor \vir (\w^3)$ with
highest weight $(\sfr{2}{3},\sfr{2}{3},\sfr{2}{3})$.
Again, the following result is easily obtained:
\begin{lem}\label{GF3A}
The Griess algebra $\GFnX{3A}$ is spanned by $\w^1$, $\w^2$, $\w^3$, $X^1$ and $X^2$.
 Moreover, we have the following commutative algebra structure on
$\GFnX{3A}$:
$$\begin{array}{c|ccccc}
{a}_{(1)}b & \w^1 & \w^2  & \w^3 & X^1  & X^2  \\ \hline
\w^1       & 2   \w^1  &  0 &  0& \frac{2}{3}X^1 &  \frac{2}{3}X^2  \\
\w^2       &    & 2   \w^2 & 0 &   \frac{2}{3}X^1 &  \frac{2}{3}X^2  \\
\w^3       &    &   & 2   \w^3 &   \frac{2}{3}X^1 &  \frac{2}{3}X^2  \\
 X^1       &    &   &      &  8X^2   & 45(\w^1+\w^2+\w^3) \\
 X^2       &    &   &      &         &        8 X^1
\end{array}
\qquad\qquad
\begin{array}{c|ccccc}
 \la a,b \ra & \w^1 & \w^2  &\w^3 &  X^1  & X^2  \\ \hline
\w^1       & \frac{2}{5}  &  0 &  0 & 0 & 0 \\
\w^2       &    & \frac{2}{5}  &  0 & 0 & 0 \\
\w^3       &    &    & \frac{2}{5}  &  0 & 0  \\
 X^1       &    &     &        &  0 & 27 \\
 X^2       &    &     &        &    &  0
\end{array}
$$
\end{lem}

\medskip

One verifies that
$$
\begin{array}{lll}
\begin{array}{l}
  \tv=\dfr{5}{14}(\w^1+\w^2+\w^3)+\dfr{1}{14} X^1
     +\dfr{1}{14} X^2,
  \vsb\\
  \tv'=\dfr{5}{14}(\w^1+\w^2+\w^3)+\dfr{\zeta}{14} X^1
     +\dfr{\zeta^{-1}}{14} X^2,
\end{array}
& &
  \la \tv,\tv'\ra = \dfr{3}{196},
  ~~~~~\mbox{where~}~\zeta =e^{2 \pi \sqrt{-1}/3}.
\end{array}
$$
In this case, the Griess algebra $\GFnX{3A}$ is not generated by
$\tv$ and $\tv'$.
Let $\nu$ be the diagram automorphism of the affine $E_6$ diagram
of order $3$ defined as:
$\alpha_0\mapsto \alpha_1 \mapsto \alpha_5\mapsto \alpha_0$,
$\alpha_6\mapsto \alpha_2 \mapsto \alpha_4 \mapsto \alpha_6$
and $\alpha_3\mapsto \alpha_3$ on the diagram \eqref{eq:4.21}.
Since $\sqrt{2}E_6$ is doubly even, we have a splitting
$\aut(V_{\sqrt{2}E_6})\simeq \hom_\Z(E_6,\C^*)\rtimes \aut(E_6)$
(see Theorem 2.1 of \cite{DN} and Chapter 5 of \cite{FLM}).
Then $\nu$ canonically acts on the Griess subalgebra above and
we find that $\tv$ and $\tv'$ generate the fixed point subalgebra
$$
  \GFnX{3A}^{\la \nu\ra}
  = \Span_\C \{ \w^1+\w^2+\w^3,X^1,X^2\} .
$$

\medskip

Summarizing, we have obtained the following table of values of
inner products between $\tv$ and $\tv'$:
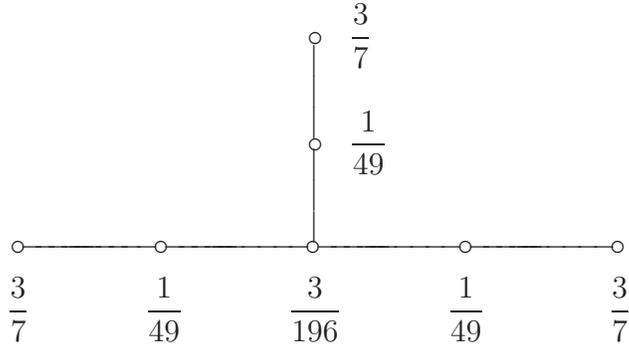
$$
\begin{array}{l}
  \hspace{120.5pt}\circ \hspace{10pt} \dfr{3}{7}
  \vspace{-15pt}\\
  \hspace{121.3pt}| \vspace{-10pt}\\
  \hspace{121.3pt}| \vspace{-10pt}\\
  \hspace{121.3pt}| \vspace{-14pt}\\
  \hspace{120.5pt}\circ \hspace{10pt}  \dfr{1}{49}
  \vspace{-16pt}\\
  \hspace{121.3pt}| \vspace{-10pt}\\
  \hspace{121.3pt}| \vspace{-10pt}\\
  \hspace{121.3pt}| \vspace{-14pt}\\
  \hspace{8pt}
  \circ\hspace{-5pt}-\hspace{-7pt}-\hspace{-7pt}-\hspace{-7pt}-

\hspace{-7pt}-\hspace{-7pt}\hspace{-7pt}-\hspace{-7pt}-\hspace{-7pt}-\hspace{
-7pt}-
  \hspace{-7pt}-\hspace{-6pt}-\hspace{-5pt}
  \circ\hspace{-5pt}-\hspace{-7pt}-\hspace{-7pt}-\hspace{-7pt}-

\hspace{-7pt}-\hspace{-7pt}\hspace{-7pt}-\hspace{-7pt}-\hspace{-6pt}-\hspace{
-6pt}-
  \hspace{-6pt}-\hspace{-6pt}-\hspace{-5pt}
  \circ\hspace{-5pt}-\hspace{-7pt}-\hspace{-7pt}-\hspace{-7pt}-

\hspace{-7pt}-\hspace{-7pt}\hspace{-7pt}-\hspace{-7pt}-\hspace{-6pt}-\hspace{
-6pt}-
  \hspace{-6pt}-\hspace{-6pt}-\hspace{-5pt}
  \circ\hspace{-5pt}-\hspace{-7pt}-\hspace{-7pt}-\hspace{-7pt}-

\hspace{-7pt}-\hspace{-7pt}\hspace{-7pt}-\hspace{-7pt}-\hspace{-6pt}-\hspace{
-6pt}-
  \hspace{-6pt}-\hspace{-6pt}-\hspace{-5pt}
  \circ
  \vspace{2pt}\\
  \hspace{7pt} \dfr{3}{7}
  \hspace{44pt} \dfr{1}{49}
  \hspace{40pt} \dfr{3}{196}
  \hspace{40pt} \dfr{1}{49}
  \hspace{47pt} \dfr{3}{7}
  \end{array}
  \vsb
$$

\begin{rem}
By the computation above, we see that in the 3A case the Griess subalgebra
generated by $\tv$ and $\tv'$ coincides with the fixed point subalgebra
$\GFnX{3A}^{\la \nu \ra}$.
This is the only case where the corresponding node is fixed by the
diagram automorphism.
\end{rem}



\section{The 3A-algebra for the Monster}\label{sec:3}

In this section, we will review and list some properties of a VOA called the
3A-algebra for the Monster which is related to certain dihedral groups of order 6 in the
Monster (cf.~\cite{LYY1,LYY2,S}). By using the VOA structure of the 3A-algebra, we will show in Theorem \ref{prop:3.16} that certain
commutant algebras of the Virasoro VOA $L(\sfr{4}{5},0)$ in an arbitrary VOA, satisfying few mild assumptions, have a subgroup of automorphisms satisfying the $3$-transposition property. These results will be
used in the last section to study the Moonshine VOA and its subalgebra related
to the Fischer group.


We first consider the extended simple Virasoro VOAs
$\W(\sfr{4}{5})=L(\sfr{4}{5},0)\oplus L(\sfr{4}{5},3)$ and
$\W(\sfr{6}{7})=L(\sfr{6}{7},0)\oplus L(\sfr{6}{7},5)$ in Theorem \ref{thm:2.7}.
It is discussed in \cite{LYY2,M3,SY} that the Moonshine VOA contains
the following subalgebra, which is a simple current extension of
$\W(\sfr{4}{5})\tensor \W(\sfr{6}{7})$:
\begin{equation}\label{eq:3.17}
\begin{array}{ll}
  U_{3A} :=& \bigl( L(\sfr{4}{5},0)\oplus L(\sfr{4}{5},3)\bigr)
 \,\tensor\, \bigl( L(\sfr{6}{7},0) \oplus L(\sfr{6}{7},5)\bigr)
 \\
  &\  \oplus\  L(\sfr{4}{5},\sfr{2}{3})^+\tensor L(\sfr{6}{7},\sfr{4}{3})^+
  \ \oplus\  L(\sfr{4}{5},\sfr{2}{3})^-\tensor L(\sfr{6}{7},\sfr{4}{3})^- .
\end{array}
\end{equation}
Moreover, a dihedral group of order $6$ can be defined using $U_{3A}$ such that
all order $3$ elements are in the 3A conjugacy class (loc. cit.). Therefore, $U_{3A}$
is closely related to the 3A-element of the Monster and we will call it the {\it 3A-algebra} for the Monster.

\begin{rem}
The 3A-algebra can be also constructed along the recipe described
in Section \ref{sec:4} via the embedding $A_2\oplus E_6\hookrightarrow E_8$,
which corresponds to the 3A node of the McKay $E_8$-observation \cite{LYY1}.
\end{rem}

In Sections \ref{sec:3.3.1} and \ref{sec:3.3.2}, we will review the results obtained
in \cite{LYY2,SY} which we will use in Section \ref{sec:3.2.3} to prove the
$3$-transposition property of Miyamoto involutions associated to derived
$c=6/7$ Virasoro vectors (cf.~Theorem \ref{prop:3.16}).


\subsection{Griess algebra}\label{sec:3.3.1}
In this subsection, we will recall some basic properties of the 3A-algebra $U_{3A}$ from \cite{LYY2,SY}.
The Griess algebra of $U_{3A}$ is $4$-dimensional and can be described as follows.

Let $\w^1$ and $\w^2$ be the Virasoro vectors of the
subalgebras $L(\sfr{4}{5},0)$ and $L(\sfr{6}{7},0)$ of
$U_{3A}$ in \eqref{eq:3.17}, respectively, and let $X^\pm$ be the
highest weight vectors of the components
$L(\sfr{4}{5},\sfr{2}{3})^\pm \tensor L(\sfr{6}{7},\sfr{4}{3})^\pm$
of $U_{3A}$.
\begin{lem}[\cite{LYY2}]\label{G3A}
The commutative algebra structure on the Griess algebra of $U_{3A}$ is given by:
$$
\begin{array}{c|cccc}
{a}_{(1)}b & \w^1 & \w^2  & X^+  & X^-  \\ \hline
\w^1       & 2   \w^1  &  0 &  \frac{2}{3}X^+ &  \frac{2}{3}X^-  \\
\w^2       &    & 2   \w^2 &   \frac{4}{3}X^+ &  \frac{4}{3}X^-  \\
 X^+       &    &          &  20X^-   & 135\w^1 + 252\w^2 \\
 X^-       &    &          &          &         20X^+
\end{array}
\qquad\qquad
\begin{array}{c|cccc}
 \la a,b \ra & \w^1 & \w^2  &  X^+  & X^-  \\ \hline
\w^1       & \frac{2}{5}  &  0 &  0 & 0 \\
\w^2       &    & \frac{3}{7}  &  0 & 0  \\
 X^+       &    &              &  0 & 81 \\
 X^-       &    &              &    &  0
\end{array}
$$
\end{lem}
The Virasoro vectors of $U_{3A}$ are classified in \cite{SY} and there are
in total three Ising vectors in $U_{3A}$.
Let $\zeta$ be a primitive cubic root of unity.
Then
\begin{equation}\label{eq:3.19}
   e^i=\dfr{5}{32} \w^1+\dfr{7}{16} \w^2 +\dfr{1}{32} \zeta^i X^+
   +\dfr{1}{32} \zeta^{-i} X^-,\q i=0,\,1,\,2,
\end{equation}
provide all the Ising vectors of $U_{3A}$. The associated $\tau$-involutions
satisfy $\abs{\tau_{e^i}\tau_{e^j}}=3$ if $i\ne j$ and therefore they generate
the symmetric group $\mathrm{S}_3$ in $\aut(U_{3A})$.
Indeed, it is known \cite{LLY,M3,SY} that $\tau_{e^i}\tau_{e^j}$ coincides
with the order three elements $\xi$ or $\xi^{-1}$ in Theorem \ref{thm:2.12}.
The VOA $U_{3A}$ is generated by any two of these Ising vectors and
correspondingly $\aut(U_{3A})\simeq \mathrm{S}_3$ is also generated by
the associated $\tau$-involutions.

It is shown in \cite{SY} that $U_{3A}$ has exactly four simple
Virasoro vectors with central charge~$4/5$, namely, $\w^1$ and the
following three vectors:
\begin{equation}\label{eq:3.20}
  x^i=\dfr{1}{16}\w^1+\dfr{7}{8}\w^2-\dfr{1}{48}\zeta^i X^+
  -\dfr{1}{48}\zeta^{-i} X^-,~~ i=0,\,1,\,2.
\end{equation}
Among these four vectors, only $\w^1$ is characteristic in the sense
it is fixed by $\aut(U_{3A})$, whereas the other three vectors are conjugated by
$\tau$-involutions $\tau_{e^i}$, $i=0$, $1$, $2$.
We call $\w^1+\w^2$ the {\it characteristic Virasoro frame\/} of $U_{3A}$.
By \eqref{eq:3.17}, we see that $\w^1$ and $\w^2$ are extendable.
Here we show that $\w^1$ is the unique extendable simple $c=4/5$ Virasoro vector
of $U_{3A}$.

\begin{lem}[\cite{LYY2}]\label{lem:3.10}
  Let $y$ be one of the $x^i$, $i=0$, $1$, $2$.
  Then as a module over $\vir(y)\simeq L(\sfr{4}{5},0)$, we have
  $U_{3A}=(U_{3A})_y[0]\oplus (U_{3A})_y[3]\oplus (U_{3A})_y[\sfr{2}{3}]
  \oplus (U_{3A})_y[\sfr{1}{8}]\oplus (U_{3A})_y[\sfr{13}{8}]$.
  Moreover, $(U_{3A})_2\cap (U_{3A})_y[\sfr{13}{8}] \ne 0$.
\end{lem}

By Theorem \ref{thm:2.12a}, there is no untwisted $\W(\sfr{4}{5})$-module
which contains an $L(\sfr{4}{5},0)$-submodule isomorphic
to $L(\sfr{4}{5},\sfr{13}{8})$, and therefore
we see that the $c=4/5$ Virasoro vectors $x^i$, $i=0$, $1$, $2$, are not extendable.


\subsection{Representation theory}\label{sec:3.3.2}

The representation theory of $U_{3A}$ was completed in \cite{SY}.

\begin{thm}[\cite{SY}]\label{thm:3.11}
  The VOA $U_{3A}$ is rational and there
  are six isomorphism types of irreducible modules over $U_{3A}$
  with the following shapes as  $\vir(\w^1)\tensor \vir(\w^2)$-modules:
  $$
  \begin{array}{lcl}
    U(0) &=& U_{3A}\ =\  [0,0]\oplus [3,0] \oplus [0,5]\oplus [3,5]\oplus
    2\, [\sfr{2}{3},\sfr{4}{3}],
    \vsb\\
    U(\sfr{1}{7}) &=& [0,\sfr{1}{7}]\oplus [0,\sfr{22}{7}]
    \oplus [3,\sfr{1}{7}] \oplus [3,\sfr{22}{7}]
    \oplus 2\, [\sfr{2}{3},\sfr{10}{21}],
    \vsb\\
    U(\sfr{5}{7}) &=& [0,\sfr{5}{7}]\oplus [0,\sfr{12}{7}]
    \oplus [3,\sfr{5}{7}] \oplus [3,\sfr{12}{7}]
    \oplus 2\, [\sfr{2}{3},\sfr{1}{21}],
    \vsb\\
    U(\sfr{2}{5}) &=& [\sfr{2}{5},0] \oplus [\sfr{2}{5},5]
    \oplus [\sfr{7}{5},0] \oplus [\sfr{7}{5},5]
    \oplus 2\, [\sfr{1}{15},\sfr{4}{3}],
    \vsb\\
    U(\sfr{19}{35}) &=& [\sfr{2}{5},\sfr{1}{7}]
    \oplus [\sfr{2}{5},\sfr{22}{7}] \oplus [\sfr{7}{5},\sfr{1}{7}]
    \oplus [\sfr{7}{5},\sfr{22}{7}] \oplus 2\, [\sfr{1}{15},\sfr{10}{21}],
    \vsb\\
    U(\sfr{4}{35}) &=& [\sfr{2}{5},\sfr{5}{7}]
    \oplus [\sfr{2}{5},\sfr{12}{7}] \oplus [\sfr{7}{5},\sfr{5}{7}]
    \oplus [\sfr{7}{5},\sfr{12}{7}] \oplus 2\, [\sfr{1}{15},\sfr{1}{21}],
  \end{array}
  $$
  where $[h_1,h_2]$ denotes an irreducible $\vir(\w^1)\tensor
  \vir(\w^2)\simeq L(\sfr{4}{5},0)\tensor L(\sfr{6}{7},0)$-module
  isomorphic to $L(\sfr{4}{5},h_1)\tensor L(\sfr{6}{7},h_2)$.
\end{thm}

By the list of irreducible modules above, we see that $U_{3A}$ is
a maximal extension of $L(\sfr{4}{5},0)\tensor L(\sfr{6}{7},0)$ as
a simple VOA.
We remark the following fundamental observation.

\begin{lem}\label{lem:3.12}
  Let $V$ be a VOA and $e\in V$ be an Ising vector.
  Then any $V$-module $M$ is $\tau_e$-stable, that is,
  the $\tau_e$-conjugated module $\tau_e \circ M $ is isomorphic
  to $M$ itself.
  In particular, if $G$ is a subgroup of $\aut(V)$ generated by
  $\tau$-involutions associated to Ising vectors of $V$ then
  $M$ is $G$-stable, that is, $g\circ M\simeq M$ for all $g\in G$.
\end{lem}

As we discussed, $\aut(U_{3A})\simeq \mathrm{S}_3$ is generated by
$\tau$-involutions associated to Ising vectors of $U_{3A}$ and
therefore all irreducible $U_{3A}$-modules are $\mathrm{S}_3$-invariant.
In general, if an irreducible $V$-module $M$ is $G$-stable then
we have a projective action of $G$ on $M$ (cf.\ \cite{DY}).
But in our case, we have an {\it ordinary\/} action of $\mathrm{S}_3$
on each irreducible $U_{3A}$-module.
For, we can find all the irreducible $U_{3A}$-modules as a submodule of
a larger VOA, say $V_{\sqrt{2}E_8}$ for example, on which we have an
ordinary $\mathrm{S}_3$-action (cf.\ \cite{LY2,LYY2}; consider
$U_{3A}^{\mathrm{S}_3}\subset (V_{\sqrt{2}A_2}\tensor V_{\sqrt{2}E_6})^+
\subset V_{\sqrt{2}E_8}^+$).
Let $M_0$, $M_1$ and $M_2$ be the principal, signature and
2-dimensional irreducible representations of $\mathrm{S}_3$.
As a $U_{3A}^{\mathrm{S}_3} \tensor \C\mathrm{S}_3$-module, one has
the following decompositions:%
\footnote{There is an Ising vector $e\in V_{\sqrt{2}E_8}$ such that
$\tau_e$ acts by $-1$ on the weight one subspace.
Then considering an embedding $U_{3A}^{\mathrm{S}_3}
\subset (V_{\sqrt{2}A_2} \tensor V_{\sqrt{2}E_6})^+\subset V_{\sqrt{2}E_8}^+$
with $e\in U_{3A}$ as in \cite{LYY2}, one can verify the decomposition.}
\begin{eqnarray}\label{eq:3.21}\nonumber
  U(0) &= &\bigl( [0,0]\oplus [3,5]\bigr) \tensor M_0
  \ \oplus\ \bigl( [3,0]\oplus [0,5]\bigr)\tensor M_1
  \ \oplus\  [\sfr{2}{3},\sfr{4}{3}]\tensor M_2
  \vsb\\ \nonumber
  U(\sfr{1}{7}) &=&
  \bigl( [0,\sfr{22}{7}] \oplus [3,\sfr{22}{7}]\bigr)\tensor M_0
  \ \oplus\ \bigl( [0,\sfr{1}{7}] \oplus [3,\sfr{1}{7}] \bigr) \tensor M_1
  \ \oplus\ [\sfr{2}{3},\sfr{10}{21}] \tensor M_2,
  \vsb\\ \nonumber
  U(\sfr{5}{7}) &=&
  \bigl( [0,\sfr{5}{7}] \oplus [3,\sfr{5}{7}] \bigr) \tensor M_0
  \ \oplus\ \bigl( [0,\sfr{12}{7}] \oplus [3,\sfr{12}{7}] \bigr)\tensor M_1
  \ \oplus\ [\sfr{2}{3},\sfr{1}{21}]\tensor M_2,
  \vsb\\ \nonumber
  U(\sfr{2}{5}) &=&
  \bigl( [\sfr{2}{5},0] \oplus [\sfr{7}{5},0]\bigr)\tensor M_0
  \ \oplus\ \bigl( [\sfr{2}{5},5]  \oplus [\sfr{7}{5},5]\bigr)\tensor M_1
  \ \oplus\ [\sfr{1}{15},\sfr{4}{3}]\tensor M_2,
  \vsb\\ \nonumber
  U(\sfr{19}{35}) &=&
  \bigl( [\sfr{2}{5},\sfr{22}{7}]  \oplus [\sfr{7}{5},\sfr{22}{7}] \bigr)\tensor M_0
  \ \oplus\ \bigl( [\sfr{2}{5},\sfr{1}{7}] \oplus [\sfr{7}{5},\sfr{1}{7}]\bigr)\tensor M_1
  \ \oplus\ [\sfr{1}{15},\sfr{10}{21}]\tensor M_2,
  \vsb\\ \nonumber
  U(\sfr{4}{35}) &=&
  \bigl( [\sfr{2}{5},\sfr{5}{7}] \oplus [\sfr{7}{5},\sfr{5}{7}] \bigr)\tensor M_0
  \ \oplus\ \bigl( [\sfr{2}{5},\sfr{12}{7}] \oplus [\sfr{7}{5},\sfr{12}{7}] \bigr)\tensor M_1
  \ \oplus\ [\sfr{1}{15},\sfr{1}{21}]\tensor M_2, \\
  & &
\end{eqnarray}
where $[h_1,h_2]$ denotes a $\vir(\w^1)\tensor \vir(\w^2)$-module
isomorphic to $L(\sfr{4}{5},h_1)\tensor L(\sfr{6}{7},h_2)$.


\subsection{$3$-transposition property of $\sigma$-involutions}\label{sec:3.2.3}

We consider involutions induced by the 3A-algebra. We will again refer to Section
\ref{sec:2.1} for the definition of simple $c=6/7$ Virasoro vectors of
$\sigma$-type and their corresponding $\sigma$-involutions.
See Lemma \ref{lem:2.5} and Eq.\ \eqref{eq:2.4} for the details.

\smallskip

Let $V$ be a VOA and let  $u\in V$ be a simple $c=4/5$ Virasoro vector.
\begin{df}\label{df:3.13}
  A simple $c=6/7$ Virasoro vector $v\in \com_V(\vir(u))$ is
  called {\it a derived Virasoro vector with respect to $u$\/} if there exists
  a sub-VOA $U$ of $V$ isomorphic to the 3A-algebra $U_{3A}$
  such that $u+v$ is the characteristic Virasoro frame of $U$.
\end{df}

\begin{lem}\label{lem:3.14}
A derived $c=6/7$ Virasoro vector $v\in \com_V(\vir(u))$ with respect to $u$ is
of $\sigma$-type on the commutant
 $\com_V(\vir(u))$.
\end{lem}

\pf
Assume that $V$ contains a subalgebra $U$
isomorphic to the 3A-algebra $U_{3A}$ as in Definition~\ref{df:3.13}.
For an irreducible $U$-module $M$, we denote
$$
  H_M:=\hom_{U_{3A}}(M,V).
$$
Then we have the isotypical decomposition
$$
  V=\bigoplus_{M\in \irr (U_{3A})} M\tensor H_M.
$$
By definition, $u+v$ is the characteristic Virasoro frame of $U_{3A}$.
Consider $V$ as a module over its subalgebra
$\vir(u)\tensor \vir(v)\tensor \com_V(U)$.
By Theorem \ref{thm:3.11} we have the following decomposition:
\begin{equation}\label{eq:3.22}
\begin{array}{l}
  V_u[0]=\vir(u)\tensor \com_V(\vir(u)),
  \vsb\\
  \com_V(\vir(u))
  = [0\oplus 5]\tensor H_{U(0)}
  \oplus [\sfr{1}{7}\oplus \sfr{22}{7}] \tensor H_{U(1/7)}
  \oplus [\sfr{5}{7}\oplus \sfr{12}{7}] \tensor H_{U(5/7)},
\end{array}
\end{equation}
where $[h_1\oplus h_2]$ denotes an irreducible $\W(\sfr{6}{7})$-module
isomorphic to $L(\sfr{6}{7},h_1)\oplus L(\sfr{6}{7},h_2)$.
By the decomposition above, we see that $v$ is of $\sigma$-type on
$\com_V(\vir(u))$.
\qed

\smallskip

We consider the one-point stabilizer
\begin{equation}\label{eq:3.23}
  \stab_{\aut(V)}(u):=\{ h\in \aut(V) \mid hu=u\} .
\end{equation}
Each $h\in \stab_{\aut(V)}(u)$ keeps the isotypical component
$V_u[h]$ invariant so that by restriction we can define a
group homomorphism
\begin{equation}\label{eq:3.24}
\begin{array}{cccc}
  \psi_u: & \stab_{\aut(V)}(u) & \longrightarrow & \aut(\com_V(\vir(u))),
  \vsb\\
  & h & \longmapsto & h|_{\com_V(\vir(u))}.
\end{array}
\end{equation}

Let $v\in \com_V(\vir(u))$ be a derived $c=6/7$ Virasoro vector with respect to $u$.
By Lemmas~\ref{lem:3.14} and \ref{lem:2.5}, we have an involution
$\sigma_v\in \aut(\com_V(\vir(u)))$.
Now let $e$ be an Ising vector of $U$.
By Lemma \ref{lem:3.12} and Eq.~\eqref{eq:3.21}, we see that $\tau_e$ keeps
$\com_V(\vir(u))$ invariant and, in fact, we have:

\begin{lem}\label{lem:3.15}
  $\psi_u(\tau_e)=\sigma_v$ for any Ising vector $e\in U$.
\end{lem}

Let $J$ be the set of all derived $c=6/7$ Virasoro vectors of $\com_V(\vir(u))$.
We will prove that the set of involutions
$$
  \{ \sigma_v \in \aut(\com_V(\vir(e))) \mid v\in J \}
$$
satisfies a $3$-transposition property.


\medskip

We say a VOA $W$ over $\R$ is {\it compact} if $W$ has a
positive definite invariant bilinear form.
A real sub-VOA $W$ is said to be a {\it compact real form} of a VOA $V$
over $\C$ if $W$ is compact and $V\simeq \C \tensor_\R W$.

We recall the following interesting theorem of Sakuma.

\begin{thm}[\cite{S}]\label{thm:3.8}
  Let $W$ be a VOA over $\R$ with grading $W=\bigoplus_{n\geq 0} W_n$,
  $W_0=\R \vac$ and $W_1=0$, and assume $W$ is compact, that is,
  the normalized invariant bilinear form on $W$ is positive definite.
  Let $x$, $y$ be Ising vectors in $W$ and denote by $U(x,y)$ the
  subalgebra of $W$ generated by $x$ and $y$. Then:
  \\
  (1) The $6$-transposition property $\abs{\tau_x\tau_y}\leq 6$ holds
  on $W$.
  \\
  (2) There are exactly nine possible inequivalent structures of
  the Griess algebra on the weight two subspace $U(x,y)_2$ of $U(x,y)$.
  \\
  (3) The Griess algebra structure on $U(x,y)_2$ is unique
  if $\abs{\tau_{x}\tau_y}=6$ and in this case $U(x,y)$ is a copy of the
  6A-algebra.
\end{thm}

\medskip
 The following is the main theorem of this section.
\begin{thm}\label{prop:3.16}
  Suppose that $V_1=0$ and $V$ has a compact real form $V_\R$ and
  every Ising vector of $V$ is in $V_\R$.
  Then for any $v^1$, $v^2\in J$, we have
  $\abs{\sigma_{v^1}\sigma_{v^2}}\leq 3$ on $\com_V(\vir(u))$.
\end{thm}

\pf Let $v^1$, $v^2 \in J$. Then there exist subalgebras $U^1$ and $U^2$ of $V$
isomorphic to the 3A-algebra such that $u$ is a simple extendable $c=4/5$
Virasoro vector in $U^1\cap U^2$ and $v^1\in \com_{U^1}(\vir(u))$, $v^2\in
\com_{U^2}(\vir(u))$. Let $a$, $a'$, $a''$ be the three distinct Ising vectors
in $U^1$ and $b$, $b'$, $b''$ be the three distinct Ising vectors in $U^2$. Set
$g:= \tau_a\tau_{a'}$. Then $g$ is an order three element induced by the
extended Virasoro VOA $\W(\sfr{4}{5})$ and we have $a'=g a$ and $a''=g^2 a$. By
our settings, $\tau_b\tau_{b'}=g$ or $g^2$ so that we may also assume
$g=\tau_b\tau_{b'}$, $b'=g b$ and $b''=g^2 b$. Note that $\psi_u(g)=1$.

Now, recall that
$$
  \sigma_{v^1}
  = \psi_u(\tau_a)
  = \psi_u(\tau_{a'})
  = \psi_u(\tau_{a''})~~~~\mbox{and}~~~~
  \sigma_{v^2}
  = \psi_u(\tau_b)
  = \psi_u(\tau_{b'})
  = \psi_u(\tau_{b''}).
$$
Thus, the order of $\sigma_{v^1}\sigma_{v^2}
=\psi_u(\tau_a \tau_{b})$ must divide that of $\tau_a \tau_b$.

By Sakuma's Theorem~\ref{thm:3.8}, we know that
$\abs{\tau_{a}\tau_{b}} \leq 6$.
If $\abs{\tau_{a}\tau_{b}}\leq 3$, there is nothing to prove.
So we assume $\abs{\tau_a \tau_b}=4$, $5$ or $6$.

First, we will note that $g$ commutes with $\tau_a \tau_b$ since $\tau_a
g=g^{-1} \tau_a$ and $\tau_b g=g^{-1} \tau_b$. \vsb

\noindent
{\bf Case $\abs{\tau_{a}\tau_{b}}=4$:}\quad
In this case, $\tau_a\tau_b g$ has order $12$, which is impossible
since $\tau_a\tau_b g = \tau_a \tau_b (\tau_b \tau_{b'})
= \tau_a \tau_{b'}$ has order $\leq 6$.
\vsb

\noindent
{\bf Case $\abs{\tau_{a}\tau_{b}}=5$:}\quad
In this case, $\tau_a\tau_b g = \tau_a \tau_{b'}$ has order $15$,
which is again impossible.
\vsb

\noindent
{\bf Case  $\abs{\tau_{a}\tau_{b}}=6$:}\quad
In this case, $\tau_a\tau_b g$ and $\tau_a\tau_b g^2$ have
order $6$ or $2$.
\vsb

\noindent
{\it Claim:}~
If $\abs{\tau_a\tau_b g}=6$ then $\abs{\tau_a\tau_b g^2}=2$.
\vsb\\
The reason is as follows.
Suppose both of them have order $6$.
Since $\tau_a\tau_bg=\tau_a\tau_{b'}$ and
$\tau_a\tau_bg^2 =\tau_a \tau_b (\tau_{b'}\tau_b)
=\tau_a \tau_{\tau_b b'}= \tau_a\tau_{b''}$, we have
$\la a,b\ra =\la a, b'\ra = \la a, b''\ra= 5/2^{10}$ and hence
$\la g^i a, g^j b\ra = 5/2^{10}$ for all $i$, $j=0$, $1$, $2$
(cf.\ the Appendix in~\cite{HLY}).

Now by using the structure of the 3A-algebra $U_{3A}$, we may
write the Ising vectors $a$ and $b$ as
$$
\begin{array}{ll}
  a &= \dfr{5}{32} u +\dfr{7}{16} v^1+ \dfr{1}{32} (X^1+X^2),
  \vsb\\
  b &= \dfr{5}{32} u +\dfr{7}{16} v^2+ \dfr{1}{32} (Y^1+Y^2),
\end{array}
$$
where $X^1$, $X^2$ (resp.\ $Y^1$, $Y^2$) are certain highest weight vectors
of weight $(\sfr{2}{3},\sfr{4}{3})$ with respect to
$\vir(u)\otimes \vir(v^1)$ (resp.\ $\vir(u)\otimes \vir(v^2)$).
Using~\eqref{eq:3.19}, we get
$$
\begin{array}{l}
  a+ga+g^2a =a+a'+a'' = \dfr{3}{32}(5 u + 14 v^1),
  \vsb\\
  b+gb+g^2b=b+b'+b''= \dfr{3}{32}(5 u + 14 v^2).
\end{array}
$$
Thus, we have
$$
  \l\la a+ga+g^2a,\, b+gb+g^2b\r\ra
  = \dfr{9}{2^{10}} \l\la 5 u + 14 v^1, 5 u  + 14 v^2\r\ra .
$$
Since $\la u, u \ra= 2/5$ and $\la u,v^i\ra =0$ for $i=1$, $2$
by Lemma~\ref{G3A}, this implies
$$
  \dfr{5}{2^{10}}
  = \frac{5^2}{2^{10}} \la u, u\ra  +\frac{7^2}{2^8} \la v^1, v^2\ra
  = \frac{5}{2^{9}} + \frac{7^2}{2^8} \la v^1, v^2\ra
$$
and thus $ \la v^1, v^2\ra= - 5/196<0$.
This is impossible by the Norton inequality
$$
  \la v^1_{(1)}v^1, v^2_{(1)}v^2  \ra
  \geq \la v^1_{(1)}v^2, v^1_{(1)}v^2 \ra \geq 0
$$
(cf.\ Theorem~6.3 and Lemma~6.5 in \cite{M1},
see also \cite{B}%
\footnote{A sketch of proof is given in \cite{B}.})
and the claim follows.

\smallskip
Therefore, $\tau_a\tau_b g=\tau_a\tau_{b'}$ or
$\tau_a\tau_b g^2=\tau_a\tau_{b''}$ has order $2$ and hence
$$
  \sigma_{v^1}\sigma_{v^2}
  = \psi_u(\tau_a \tau_{b})
  = \psi_u(\tau_a \tau_{b'})
  = \psi_u(\tau_a \tau_{b''})
$$
is of order at most~$2$.
\qed


\section{The Fischer group}\label{sec:5}

In this section, we will discuss the properties of the commutant vertex operator
subalgebra $\VF$ of the Moonshine VOA $V^\natural$. We will show that the
full Fischer $3$-transposition group $\Fi$ is a subgroup of $\aut(\VF)$. We also
show that there exist one-to-one correspondences  between 2C-involutions of
$\Fi$ and derived $c=6/7$ Virasoro vectors in $\VF$.

Finally, we will discuss the embeddings of $\UFnX{nX}$ into $\VF$. The main idea
is to embed the root lattice $E_6$ into $E_8$ and view the $\UFnX{nX}$ as
certain commutant subalgebras of  the lattice VOA $V_{\sqrt{2}E_8}$. Then we
shall show that the product of two $\sigma$-involutions generated by $c=6/7$
Virasoro vectors in $\UFnX{nX}$ exactly belong to the conjugacy class $nX$ in
$\Fi$.
By this procedure, we obtain a VOA description of the $E_6$
structure inside ${\rm Fi}_{24}$.

\smallskip

The automorphism group of the Moonshine VOA $V^\natural$ is the Monster
simple group $\M$  \cite{FLM}. Consider the monstrous Griess algebra of
dimension $196884$~\cite{C,G}. It is known that the monstrous Griess algebra is
naturally realized as the subspace of weight~$2$ of $V^\natural$ \cite{FLM},
which we call the Griess algebra of $V^\natural$ and denote by $\BN$. We will
freely use the character tables in \cite{ATLAS}.



\subsection{The Fischer group vertex operator algebra $\VF$}\label{sec:5.2.1}

We denote by $\mathrm{Fi}_{24}$ the Fischer $3$-transposition
group and by $\mathrm{Fi}_{24}'$ its derived subgroup,
the 3rd largest sporadic finite simple group.
Let $g\in \M$ be a 3A-element.
Then $C_\M(g)\simeq 3.\mathrm{Fi}_{24}'$ and it is shown in
\cite{C,MN} that the monstrous Griess algebra $\BN$ has
an irreducible decomposition
\begin{equation}\label{eq:5.9}
  \BN
  =\ \ul{\bf 1}\ \oplus\ \ul{\bf 1}\ \oplus\ \ul{\bf 8671}\ \oplus\
  \ul{\bf 57477}\ \oplus\ 2\cd \ul{\bf 783}\ \oplus\ 2\cd \ul{\bf 64584}\
\end{equation}
as a $C_\M(g)$-module.
Therefore, we can take a Virasoro vector $u\in (\BN)^{C_\M(g)}$
such that $(\BN)^{C_\M(g)}=\C u \oplus \C (\w -u)$ is an
orthogonal sum.
We take $u$ to be the shorter one, that is, the central charge
of $u$ is smaller than that of $\w-u$.
It is also shown in loc.\ cit.\ that the central charge of
the shorter Virasoro vector $u$ is $4/5$ and its spectrum on $\BN$
is as follows:
\begin{equation}\label{eq:5.10}
\begin{array}{ccccccccccccc}
  \BN
  &=& \ul{\bf 1} &\oplus& \ul{\bf 1} &\oplus& \ul{\bf 8671} &\oplus&
  \ul{\bf 57477} &\oplus& 2\cd \ul{\bf 783} &\oplus& 2\cd \ul{\bf 64584}
  \vsb\\
  u_{(0)} & : & 0 && 2 && 2/5 && 0 && 2/3 && 1/15
\end{array}
\end{equation}
The following result about the extendibility seems already known
to experts (see for example \cite{MN,KMY,M2}), even though no rigorous proof
has been given so far.

\begin{thm}\label{thm:5.16}
  Let $g$ be a 3A-element of the Monster.
  The cyclic group $\la g\ra$ uniquely determines an extendable
  simple $c=4/5$ Virasoro vector in $(V^\natural)^{C_\M(g)}$.
\end{thm}

\pf
By the decomposition in Eq.~\eqref{eq:5.9}, every cyclic subgroup $\la g\ra$ defines a unique simple
$c=4/5$ Virasoro vector in $(V^\natural)^{C_\M(g)}$. We will prove that this
Virasoro vector is extendable.

Let $t$ be a 2A-involution in $N_\M(\la g\ra )$ ($\simeq 3.\Fi$) but not in
$C_\M(g)$ ($\simeq 3.\Fi'$).
Then the subgroup $H$ generated by $t$ and $g$ in $\M$ is isomorphic to
$\mathrm{S}_3$ and $C_\M(H)\simeq \mathrm{Fi}_{23}$ (cf. \cite{ATLAS} and Lemma 13.3 of \cite{G}).
Take an involution ${t'}\in H$ such that $t{t'}=g$.
Then ${t'}$ is conjugate to $t$ and $H$ is generated by $t$ and ${t'}$.
By the one-to-one correspondence between 2A-elements of $\M$ and
Ising vectors of $V^\natural$
(cf.~\cite{M1} and \cite{Ho}, Lemma~3; see also \cite{HLY}, Theorem.~5.1),
there exist Ising vectors
$e^0$, ${e^1} \in V^\natural$ such that $\tau_{e^0}=t$,
$\tau_{{e^1}}= {t'}$ and $e^0$ and ${e^1}$ are fixed by
$C_\M(t)$ and $C_\M({t'})$, respectively.
Since $C_\M(H)$ is a subgroup of both $C_\M(t)$ and $C_\M({t'})$,
$e^0$ and ${e^1}$ are both contained in $(V^\natural)^{C_\M(H)}$.
By \cite{ATLAS}, one obtains the following decomposition of
the Griess algebra $\BN$ as a $C_\M(H)$-module:
\begin{equation}\label{eq:5.11}
  \BN
  = 5 \cd \ul{\bf 1} \oplus 3 \cd \ul{\bf 782} \oplus
    3 \cd \ul{\bf 3588} \oplus \ul{\bf 5083} \oplus \ul{\bf 25806}
    \oplus \ul{\bf 30888} \oplus 2 \cd \ul{\bf 60996}.
\end{equation}
Let $U$ be the subalgebra of $V^\natural$ generated by
$e^0$ and ${e^1}$.
Then $U$ is isomorphic to the 3A-algebra $U_{3A}$ \cite{M3,SY}.
Since the Griess algebra of $U$ is 4-dimensional
(cf.\ Section~\ref{sec:3.3.1}), it follows from the decomposition above
that the weight two subspace of $(V^\natural)^{C_\M(H)}$ is spanned
by that of $U$ and the conformal vector of $V^\natural$.
Now it is clear that all simple $c=4/5$ Virasoro vectors of
$(V^\natural)^{C_\M(H)}$ are contained in $U$.
The 3A-algebra has one extendable simple $c=4/5$ Virasoro vector
and three non-extendable ones.
By Lemma~\ref{lem:3.10}, the non-extendable ones in the 3A-algebra
have eigenvalues $13/8$ on the Griess algebra of the 3A-subalgebra
and these eigenvalues do not appear in the decomposition~\eqref{eq:5.10}.
Therefore, the subalgebra $(V^\natural)^{C_\M(g)}$ contains the unique
simple $c=4/5$ Virasoro vector, which is extendable as claimed.
\qed
\vsb

Let $e^0$ and $e^1$ be Ising vectors as in the proof above.
Without loss, we may assume that the 3A-element $\xi_u$ associated to $u$,
defined as in Theorem~\ref{thm:2.12}, coincides with $\tau_{e^0}\tau_{e^1}$.

It follows from Theorem \ref{thm:2.12} (cf.\ \cite{KMY,M2}) that for each
embedding $\W(\sfr{4}{5})= L(\sfr{4}{5},0)\oplus L(\sfr{4}{5},3) \hookrightarrow
V^\natural$, one can define an element $\xi$ of the Monster with $\xi^3=1$. By
computing its trace on $V^\natural_2$, one can show that $\xi$ belongs to the
conjugacy class 3A (cf.~Section 4.1 of \cite{Ma1}). We remark here that a single
$L(\sfr{4}{5},0)$ does not define an order three symmetry, and in order to obtain
3A-elements, we have to extend $L(\sfr{4}{5},0)$ to a larger algebra
$\W(\sfr{4}{5})=L(\sfr{4}{5},0)\oplus L(\sfr{4}{5},3)$ (cf.\ \cite{M2}).
It remains a natural question whether the map associating the order three
element $\xi$ defined as in Theorem \ref{thm:2.12} is injective or not.
Since the order three element is defined not only by the Griess algebra
but also by the $3$-primary vector, we would need extra information about
the weight three subspace of $V^\natural$ to solve this question.


Because of this problem, we first fix a 3A-element $g\in \M$ and
then take the unique simple extendable $c=4/5$ Virasoro vector $u$
in the fixed point subalgebra $(V^\natural)^{C_\M(g)}$. Then
$g=\xi_u$ by the preceding argument. Let $w$ be an associated
$3$-primary vector of $u$ and denote by $\W(u,w)$ the sub-VOA
generated by $u$ and $w$ which is isomorphic to the extended
Virasoro VOA $\W(\sfr{4}{5})$.

\begin{df}\rm
The {\it Fischer group vertex operator algebra\/} is defined as the commutant
\begin{equation}\label{eq:5.12}
  \VF :=\com_{V^\natural}(\W(u,w))= \com_{V^\natural}(\vir(u)) .
\end{equation}
\end{df}

Below we will show that the VOA $\VF$ affords an action of the Fischer
$3$-transposition group $\Fi$.

\begin{lem}\label{lem:5.17}
  $N_\M(\la \xi_u\ra )\subset \stab_\M(u)$.
\end{lem}

\pf
Since $\aut(\la\xi_u\ra)\simeq \aut(\Z_3)=\Z_2$, we have either
$N_\M(\la \xi_u\ra )=C_\M(\xi_u)$ or $C_\M(\xi_u)$ is normal of index $2$
in $N_\M(\la \xi_u\ra )$.
We have seen in Theorem \ref{thm:5.16} that $C_\M(\xi_u)\subset \stab_\M(u)$.
For elements $h\in C_\M(\xi_u)$ and $g \in N_\M(\la \xi_u\ra )$
we get $hgu =gh^{g}u=gu$ where $h^g= g^{-1}hg \in C_\M(\xi_u)$.
Thus $gu$ is fixed by $C_\M(\xi_u)$. By the decomposition \eqref{eq:5.9}
of $\BN$ as a $C_\M(\xi_u)$-module we see that $gu$ must be
equal to the shorter Virasoro element $u$ in the fixed point subspace
${(\BN)}^{C_\M(\xi_u)}$.
Thus $N_\M(\la \xi_u\ra )\subset \stab_\M(u)$.
\qed
\vsb

Recall the group homomorphism $\psi_u$ defined in \eqref{eq:3.24}.
By Lemma \ref{lem:5.17}, we have an action of $N_\M(\la \xi_u\ra )$
on $\com_{V^\natural}(\vir(u))$ via $\psi_u$.

\begin{prop}\label{prop:5.18}
  Consider the homomorphism
  $\psi_u : \stab_\M(u) \to \aut(\com_{V^\natural}(\vir(u)))$.
  Then $\psi_u(N_\M(\la \xi_u\ra ))\simeq \mathrm{Fi}_{24}$.
\end{prop}
\pf
It is clear that $\la \xi_u\ra\subset \ker \psi_u$.
By \eqref{eq:5.10}, it is also clear that $\psi_u(C_\M(\xi_u))\ne 1$.
Therefore, either $\psi_u(N_\M(\la \xi_u\ra ))\simeq \mathrm{Fi}_{24}'$
or $\mathrm{Fi}_{24}$.
Since $\tau_{e^0}\in N_\M(\la \xi_u\ra )\setminus C_\M(\xi_u)$ acts
non-trivially on $\com_{V^\natural}(\vir(u))$ via $\psi_u$,
we see $\psi_u(N_\M(\la \xi_u\ra ))\simeq \mathrm{Fi}_{24}$.
\qed

\smallskip

By \eqref{eq:5.9} and \eqref{eq:5.10} the Griess algebra $\VF_2$
is of dimension $57478$ and has, under $\psi_u(N_\M(\la \xi_u\ra ))\simeq \Fi$,
the decomposition
\begin{equation}\label{eq:5.13}
  \VF_2 = \ul{\bf 1}\oplus \ul{\bf 57477}.
\end{equation}

\medskip

We have seen in Proposition \ref{prop:5.18} that $\VF$ naturally affords an
action of the full Fischer group $\Fi$. We cannot prove that $\Fi$ is the full
automorphism group of $\VF$ at the moment. Instead, we will prove the
following partial result.

\begin{thm}\label{thm:5.19}
  Let $\EuScript{X}$ be the subalgebra of $\VF$ generated by
  the weight $2$ subspace.
  Then $\aut(\EuScript{X})\simeq N_\M(\la \xi_u\ra )/\la \xi_u\ra
  \simeq \mathrm{Fi}_{24}$.
\end{thm}

\pf
It is shown in \cite{GL, LM} that there exists a $14$-dimensional
sublattice $L$ of the Leech lattice  $\Lambda$ such that $V_L^+$ contains
a sub-VOA isomorphic to $U_{3A}$.
It is also known that the annihilator
$\mathrm{Ann}_\Lambda(L)=\{ \alpha\in \Lambda \mid \la \al, \beta\ra =0
\text{ for all }\beta \in L\}$ contains a sublattice isomorphic
to $\sqrt{2}A_2^{\oplus 5}$ (cf.~\cite{GL}).
Therefore we can find a sub-VOA  $U_{3A}\tensor V_{\sqrt{2}A_2}^+$
of $V_\Lambda^+$, which is also contained in the Moonshine VOA.
This shows that $\VF$ contains a sub-VOA isomorphic to
$V_{\sqrt{2}A_2}^+$.
It is well-known that $V_{\sqrt{2}A_2}^+$ contains 6 Ising vectors
(cf.\ \cite{LSY}) and therefore $\VF$ contains at least 6 Ising vectors.

We have already seen that $\psi_u(N_\M(\la \xi_u\ra ))\simeq \mathrm{Fi}_{24}$
faithfully acts on $\EuScript{X}$.
Let $e^0$ and $e^1$ be Ising vectors of $V^\natural$ such that
$\tau_{e^0}\tau_{e^1}=\xi_u$.
Then $e^0$ and $e^1$ generate a subalgebra $U(e^0,e^1)$ isomorphic
to the 3A-algebra such that $\W(u,w)\subset U(e^0,e^1)$.
We can take an Ising vector $e^2$ of $V^\natural$
orthogonal to both $e^0$ and $e^1$.
Then $g e^2\in \EuScript{X}$ and $\tau_{ge^2} \in C_\M(\xi_u)$ for any
$g\in \aut(\EuScript{X})$.
As we have shown, $\ker \psi_u = \la \xi_u \ra$ and hence
$\{ \psi_u(\tau_{ge^2})=g\psi_u(\tau_{e^2}) g^{-1} \mid g\in
\aut(\EuScript{X})\}$ define non-trivial involutions on $\EuScript{X}$.
Thus the subgroup generated by
$\{ \psi_u(\tau_{g e^2}) \mid g\in \aut(\EuScript{X})\}$ is normal
in $\aut(\EuScript{X})$ and isomorphic to
$C_\M(\xi_u)/\la \xi_u\ra\simeq \mathrm{Fi}_{24}'$.

Now by conjugation we can define a group homomorphism
$\alpha: \aut(\EuScript{X}) \to \aut (\mathrm{Fi}_{24}') \simeq
\aut(\psi_u(C_\M(\xi_u)))$.
Since $\mathrm{Out}(\mathrm{Fi}_{24}')=2$ by \cite{ATLAS} and
$\psi_u(N_\M(\la \xi_u\ra ))$ contains an outer involution defined by
a simple $c=6/7$ Virasoro vector $v\in U(e^0,e^1)\cap \EuScript{X}$,
we see that
$$
  \alpha(\aut(\EuScript{X}))
  = \alpha(\psi_u(N_\M(\la\xi_u\ra )))
  = \aut(\mathrm{Fi}_{24}')
  \simeq \mathrm{Fi}_{24} .
$$
This implies $\aut(\EuScript{X})\simeq \ker \alpha \times
\mathrm{Fi}_{24}$.
The Griess algebra of $\VF$ is $57478$-dimensional
and it has a decomposition $57478=\ul{\bf 1}\oplus \ul{\bf 57477}$ as
a module over $\psi_u(C_\M(\xi_u))\simeq \mathrm{Fi}_{24}'$.
So the Griess algebra of $\VF$ is spanned by its conformal vector
$\w_{\VF}$ and Ising vectors $\{ g e^2 \mid g\in C_\M(\xi_u)\}$.
Now let $h\in \ker\alpha$.
Then $h$ acts by a scalar on the 57477-dimensional component, say
$\lambda$.
Write $e^2=p\w_{\VF} +x$ with $p\in \C$ and $x\in \ul{\bf 57477}$.
Then $x\ne 0$ since the central charge of $\w_{\VF}$ is equal to
$24-4/5$ and $he^2=p\w +\lambda x$.
Since both $e^2/2$ and $h e^2 /2$ are idempotents in
the Griess  algebra, it follows that $\lambda=1$ and $h=1$.
Therefore $\ker \alpha =1$ and we obtain the desired isomorphism
$\aut(\EuScript{X})=\psi_u(N_\M(\la \xi_u\ra ))\simeq \mathrm{Fi}_{24}$.
\qed


\subsection{The $3$-transposition property}

In this section, we will establish the correspondence between derived $c=6/7$
Virasoro vectors in $\VF$ and 2C-involutions of the Fischer group $\Fi$.

Take a 2C-involution $t$ of $\Fi$ and consider the decomposition of
$\VF$ as a module over $C_{\Fi}(t)\simeq \mathrm{Fi}_{23}$.
In the computation of \eqref{eq:5.11} we have already obtained that
\begin{equation}\label{eq:5.14}
  \VF_2 = \ul{\bf 1} \oplus \ul{\bf 1} \oplus \ul{\bf 30888} \oplus \ul{\bf 25806}
  \oplus \ul{\bf 782}
\end{equation}
as a $C_{\Fi}(t)$-module.
Therefore, the fixed point subalgebra $(\VF_2)^{C_{\Fi}(t)}$ is
2-dimensional and forms a commutative associative algebra spanned by two
mutually orthogonal Virasoro vectors.
In order to determine the central charge of the shorter Virasoro
vector in $(\VF)^{C_{\Fi}(t)}$, we use the 3A-algebra for the Monster
to obtain suitable decompositions.

Let $v$ be a derived $c=6/7$ Virasoro vector in $\com_{V^\natural}(\vir(u))$
with respect to $u$. Let $U\subset V^\natural$ be the corresponding sub-VOA
isomorphic to $U_{3A}$, that is, $u+v$ is the conformal vector of $U$, and let
$e^0$, $e^1$ and $e^2$ be Ising vectors of $U$ such that
$\tau_{e^0}\tau_{e^1}=\xi_u$.

\smallskip

For an irreducible $U$-module $M$, we set
$H_M:=\hom_{U}(M,V^\natural)$.
Then we have the decomposition
\begin{equation}\label{eq:5.15}
  V^\natural = \bigoplus_{M\in \irr(U_{3A})} M\tensor H_M.
\end{equation}
Clearly $H_M$ forms a module over the commutant subalgebra
$\com_{V^\natural}(U)$.
\begin{lem}
The top weight $h(H_M)$  and the dimension $d(H_M)$ of the top level of the
$\com_{V^\natural}(U)$-modules $H_M$ are given by the following table
\begin{equation}\label{eq:5.16}
\begin{array}{|c||c|c|c|c|c|c|}\hline
  M & U(0) & U(\sfr{1}{7}) & U(\sfr{5}{7}) & U(\sfr{2}{5})
  & U(\sfr{19}{35}) & U(\sfr{4}{35})
  \\ \hline
  h(H_M) & 0 & 13/7 & 9/7 & 8/5 & 51/35 & 66/35
  \\ \hline
  d(H_M) & 1 & 25806 & 782 & 5083 & 3588 & 60996
  \\ \hline
\end{array}
\end{equation}
Moreover, one has $\dim (H_{U(0)})_2=30889$.
\end{lem}
\pf
By Lemma \ref{lem:2.8} and the classification of irreducible
$U_{3A}$-modules in Theorem~\ref{thm:3.11}, we know that the possible
eigenvalues of $u_{(1)}$ on $\BN$ are $0$, $\sfr{1}{15}$, $\sfr{2}{5}$ and $\sfr{2}{3}$,
and those for $v_{(1)}$ are $0$, $\sfr{1}{7}$, $\sfr{5}{7}$, $\sfr{4}{3}$, $\sfr{1}{21}$ and $\sfr{10}{21}$.
Applying a similar computation as in the proof of Lemma~5.2 of~\cite{HLY},
we obtain the Lemma.
\qed
\vsb

On $\com_{V^\natural}(\vir(u))$, one can define  the $\sigma$-involution
$\sigma_v$ as in \eqref{eq:2.4}, which coincides with
$\psi_u(\tau_{e^0})$ by Lemma \ref{lem:3.15}.

\begin{prop}\label{prop:5.21}
 The involution
 $\sigma_v=\psi_u(\tau_{e^0})$ is a 2C-element of
  $\psi_u(N_\M(\la \xi_u\ra ))\simeq \mathrm{Fi}_{24}$.
\end{prop}

\pf
Let us consider the trace of $\sigma_v=\psi_u(\tau_{e^0})$
on the Griess algebra of $\VF =\com_{V^\natural}(\vir(u))$.
By \eqref{eq:5.16} and \eqref{eq:3.22}, one has
$$
  \mathrm{Tr}_{\VF_2}\sigma_v
  = 1+\dim (H_{U(0)})_2 +\dim H_{U(5/7)}-\dim H_{U(1/7)} = 5866,
$$
which coincides only with the trace of a 2C-involution of $\mathrm{Fi}_{24}$ on
$\VF_2=\ul{\bf 1} \oplus \ul{\bf 57477}$ by \cite{ATLAS}.
\qed

\begin{lem}\label{lem:5.22}
  The set $\{e^0, e^1,e^2\}$ consisting of the three Ising vectors of $U$
  is stabilized by
  $\psi_u^{-1}(C_{\Fi}(\psi_u(\tau_{e^0})))$.
\end{lem}

\pf
Take any $g\in \psi_u^{-1}C_{\Fi}(\psi_u(\tau_{e^0})) \subset N_\M(\la \xi_u\ra )$.
Then $[g, \tau_{e^0}]\in \la \xi_u \ra$ and  hence
$\tau_{ge^0}=g\tau_{e^0} g^{-1} \in \{ \tau_{e^0}, \xi_u \tau_{e^0},
 \xi_u^2 \tau_{e^0}\}$.
Since $\tau_{e^0}\tau_{e^1}= \xi_u$ we get
$\{\xi_u\tau_{e^0}, \xi_u^2\tau_{e^0}\}=\{\tau_{e^1}, \tau_{e^2}\}$.
Thus, $ge^0\in \{e^0, e^1,e^2\}$ by the one-to-one
correspondence~\cite{M1,Ho}. Similarly, we also have
$ge^1\in \{e^0, e^1, e^2\} $ and $ge^2 \in \{e^0, e^1, e^2\}$.
\qed

\begin{prop}\label{prop:5.23}
  A derived $c=6/7$ Virasoro vector $v\in \com_{V^\natural}(\vir(u))$
  with respect to $u$ is fixed by the centralizer of
  the 2C-involution $\psi_u(\tau_{e})$ of the Fischer group $\Fi$.
\end{prop}

\pf Let $\al =e^0+e^1+e^2$. Then, $\al$ is fixed by
$\psi_u^{-1}C_{\Fi}(\psi_u(\tau_{e^0}))$ by Lemma \ref{lem:5.22}. On the other
hand, by Lemma~\ref{G3A} and \eqref{eq:3.19}, we have
\[
\al =e^0+e^1+e^2 =\dfr{15}{32}u +\dfr{21}{16}v \quad \text{ and }\quad
\al^2= 2\left(\dfr{15}{32}\right)^2 u +2 \left(\dfr{21}{16}\right)^2 v,
\]
where $v$ is the derived $c=6/7$ Virasoro vector in $U$. Thus,
$$
  v=  \frac{16}{567}(16 \alpha^2 -15 \alpha).
$$
Hence, $v$ is fixed by the centralizer of the 2C-involution
$\psi_u(\tau_{e^0})$ in $\Fi$.
\qed
\smallskip

Now we establish the one-to-one correspondence between
2C-involutions of $\Fi$ and derived $c=6/7$ Virasoro
vectors of $\VF$, which is one of our main results.

\begin{thm}\label{thm:5.25}
  The map which associates  a derived  $c=6/7$ Virasoro vector
  to its \hbox{$\sigma$-involution} defines a bijection between
  the set of all derived  $c=6/7$ Virasoro vectors in
  $\com_{V^\natural}(\vir(u))$ with respect to $u$ and
  the 2C-conjugacy class of $\Fi=\psi_u(N_\M(\la \xi_u\ra ))$.
\end{thm}

\pf The map of the theorem is equivariant with respect to the natural action of
$\psi_u(N_\M(\la \xi_u\ra ))$ on the derived  vectors and the conjugation action of
$\Fi$ on the set of its 2C-involutions, respectively.

As seen in the proof of Theorem~\ref{thm:5.16}, the vector $u$ is contained in
a sub-VOA isomorphic to the 3A-algebra $U_{3A}$. Thus there exists at least one
derived  $c=6/7$ Virasoro vector with respect to $u$ in
$\com_{V^\natural}(\vir(u))$. The transitivity of the conjugation action on the
2C-involutions shows now the surjectivity of the map.

For the injectivity, fix a 2C-involution $t$ of $\Fi=\psi_u(N_\M(\la \xi_u\ra ))$. By
Proposition~\ref{prop:5.23}, any derived $c=6/7$ Virasoro vector $v$ in
$\com_{V^\natural}(\vir(u))$ such that $\sigma_v= t$ is contained in $
(\VF)^{C_{\Fi}(t)}$. We have seen in \eqref{eq:5.14} that the fixed point
subalgebra $(\VF_2)^{C_{\Fi}(t)}$ is spanned by two mutually orthogonal
Virasoro vectors. Hence $v$ must be the unique shorter Virasoro vector of
$(\VF)^{C_{\Fi}(t)}$ of central charge $c=6/7$. \qed

The proof gives also the following corollary.
\begin{cor}\label{cor:5.24}
  Every 2C-involution $t$ of the Fischer group $\Fi$ defines an
  unique derived $c=6/7$ Virasoro vector of the fixed
  point subalgebra $(\VF)^{C_{\Fi}(t)}$.
\end{cor}

As a consequence of Theorem \ref{prop:3.16} and
Theorem~\ref{thm:5.25}, we also recover:
\begin{cor}\label{cor:5.25}
The $2C$ involutions of the Fischer group $\Fi$ satisfy the $3$-transposition
property.
\end{cor}

\begin{rem}\label{rem:5.27}
  We expect that the full automorphism group of $\VF$ is actually $\Fi$.
  Because of Theorem \ref{thm:5.19}, it suffices to show
  $\VF = \EuScript{X}$, that is, $\VF$ is generated by its Griess algebra.
  This is a technically difficult problem since we do not know a nice
  embedding of $\W(\sfr{4}{5})$ into $V^\natural$ to study the commutant
  subalgebra $\VF = \com_{V^\natural}(\W(\sfr{4}{5}))$.
  It is conjectured in \cite{DM} that one can obtain $V^\natural$ from $V_\Lambda$
  by the $\Z_3$-orbifold construction where the $\Z_3$-automorphism is
  induced by a
 automorphism of the $A_2$ lattice via
  an embedding $\sqrt{2}A_2^{12}\hookrightarrow \Lambda$.
  If this conjectural $\Z_3$-orbifold construction is established, we obtain
  a natural embedding of $\W(\sfr{4}{5})$ into $V^\natural$ and then
  we can solve the question above immediately by a decomposition of $V^\natural$
  given in \cite{KLY}.
\end{rem}


\subsection{Embedding of $\UFnX{nX}$ into $\VF$}\label{sec:5.2.3}


We recall the definition of $\VFnX{nX}$ from Section~\ref{sec:4.1}.
We have a full rank sublattice
$Q\oplus E_6\simeq A_2\oplus E_6$ of $E_8$.
Since the index of $A_2\oplus E_6$ in $E_8$ is three,
we have a coset decomposition
$$
  E_8= A_2\oplus E_6 \sqcup (\delta +A_2\oplus E_6) \sqcup
  (2\delta +A_2\oplus E_6)
$$
with some $\delta \in E_8$ and correspondingly we obtain a decomposition
$$
  V_{\sqrt{2}E_8}
  = V_{\sqrt{2}(A_2\oplus E_6)} \oplus V_{\sqrt{2}(\delta+ A_2\oplus E_6)}
  \oplus V_{\sqrt{2}(2\delta +A_2\oplus E_6)}.
$$
Define $\eta \in \aut(V_{\sqrt{2}E_8})$ by
$$
  \eta =
  \begin{cases}
  1 & ~~\mbox{on}~~~V_{\sqrt{2}(A_2\oplus E_6)},
  \\
  e^{2\pii/3} & ~~\mbox{on}~~~V_{\sqrt{2}(\delta +A_2\oplus E_6)},
  \\
  e^{4\pii/3} & ~~\mbox{on}~~~V_{\sqrt{2}(2\delta +A_2\oplus E_6)}.
  \end{cases}
$$
Then $\eta$ is clearly in $\mathcal{F}_{nX}$, see \eqref{eq:4.9}.
Indeed, $\mathcal{F}_{nX}$ is generated by $\eta$ and $\rho_{nX}$.
Note that we can write down $\rho_{nX}$ in exponential form
$$
  \rho_{nX}=\exp(2\pii \gamma^{nX}_{(0)}/n)\quad \text{ with suitable }
  \gamma^{nX} \in L_{nX}^*,
$$
which also defines an automorphism of $V_{\sqrt{2}E_8}$ and fixes
$V_{\sqrt{2}\EL_{nX}}$ pointwisely.

\begin{rem}\label{rem:5.28}
  Recall $\tilde{\om}_{Q}$ with $Q \simeq A_2$ is a simple extendable $c=4/5$
  Virasoro vector in a lattice VOA $V_{\sqrt{2}Q}$ and $\UFnX{nX}$ equals
  the commutant subalgebra $\com_{\VFnX{nX}}(\vir(\tilde{\om}_{Q}))$
  in $\VFnX{nX}$.
  Moreover, $\rho_{nX}$ fixes $\tilde{\om}_{Q}$.
  Let $U^1$ be the subalgebra generated by $\hat{e}$ and $\eta(\hat{e})$,
  and set $U^2=\rho_{nX}( U^1)$.
  Then $U^1\simeq U^2\simeq U_{3A}$ and $\tilde{\om}_{Q}$ is contained in both
  $U^1$ and $U^2$.
  Note that $\tv=\tilde{\w}_{E_6}$ is contained in
  $\com_{U^1}( \vir(\tilde{\om}_{Q}))$.
  Similarly, $\tv'= \rho_{nX}\tv$ of $\UFnX{nX}$ is contained in
  $\com_{U^2}( \rho_{nX}(\vir(\tilde{\om}_{Q}))) =
  \com_{U^2}( \vir(\tilde{\om}_{Q}))$.
  Thus, the $c=6/7$ Virasoro vectors $\tv$ and   $\tv'$ of $\UFnX{nX}$ are
  derived Virasoro vectors with respect to $\tilde{\om}_{Q}$.
\end{rem}

\begin{prop}\label{prop:5.29}
  For any $nX=1A$, $2A$ or $3A$, the VOA $\VFnX{nX}$ can be embedded into
  the Moonshine VOA $V^\natural$.
\end{prop}

\pf
As we have shown in Section \ref{sec:3.3new},
$\VFnX{1A}$ is isomorphic to the monstrous 3A-algebra $U_{3A}$ and
$\VFnX{2A}$ is isomorphic to the monstrous 6A-algebra $U_{6A}$
discussed in \cite{LYY2}.
It is shown in \cite{LM} that both $U_{3A}$ and $U_{6A}$ are subalgebras
of $V^\natural$ and therefore $\VFnX{1A}$ and $\VFnX{2A}$ are also
contained in $V^\natural$.
That $\VFnX{3A}\simeq M_{\mathcal{C}_4}$ is contained in $V^\natural$
will be shown in Appendix~\ref{sec:A.4}.
\qed

\medskip

Finally, we will establish our main theorem.

\begin{thm}\label{thm:5.30}
  Let $u$ be a simple extendable $c=4/5$ Virasoro vector in $V^\natural$
  such that $u\in (V^\natural)^{C_\M(\xi_u)}$.
  Then for any $nX=1A$, $2A$, $3A$, the VOA $\UFnX{nX}$
  can be embedded into $\VF= \com_{V^\natural}(\vir(u))$.
  Moreover, $\sigma_{\tv}\sigma_{\tv'}$ belongs to the
  conjugacy class $nX$ of $\Fi \subset \aut(\VF)$.
\end{thm}

\pf
First we embed $\VFnX{nX}$ into $V^\natural$ using Proposition \ref{prop:5.29}.
Let $e$, $e'$  be a pair of Ising vectors of $\VFnX{nX}$ that generate
a subalgebra $U$ such that $\tilde{\om}_{Q}\in U$ and $U\simeq U_{3A}$.
Since pairs of Ising vectors in $V^\natural$ generating the 3A-algebra are
mutually conjugate under $\aut(V^\natural)$, we may identify
$\tilde{\om}_{Q}$ with $u$ by Theorem~\ref{thm:5.16}.
Thus, we have
$$
  \UFnX{nX}
  = \com_{\VFnX{nX}}(\vir(\tilde{\om}_{Q(E_6)}))\subset \com_{V^\natural}(\vir(u))
  = \VF
$$
as desired.

Next we will show that $h:=\sigma_{\tv}\sigma_{\tv'}$ belongs to the class $nX$
of $\mathrm{Fi}_{24}= \aut(\EuScript{X})$. Note that $\tv$, $\tv'\in \VF_2
\subset \EuScript{X}$. Recall that there is an exact sequence
$$
  1 \longrightarrow \la \xi_u\ra \longrightarrow N_\M(\la \xi_u\ra )
  \longrightarrow \aut(\EuScript{X})\simeq \mathrm{Fi}_{24}
  \longrightarrow  1
$$
with the projection map $\psi_u: N_\M(\la \xi_u\ra ) \to \aut(\EuScript{X})$.
Let $e^1$ and $e^2$ be Ising vectors in $V_{nX}$ such that
$\psi_u(\tau_{e^1})= \sigma_{\tv}$ and
$\psi_u(\tau_{e^2})= \sigma_{\tv'}$.
Set $g=\tau_{e^1}\tau_{e^2}$.
Then $h=\psi_u(g)$ and the inverse image
$\psi_u^{-1}( \la h \ra)$ has order $3n$ and is generated by
$\xi_u$ and $g$.

\smallskip
\noindent
{\bf 1A case:}  In this case, $\tv=\tv'$ and
hence $h=\sigma_{\tv}\sigma_{\tv'}$ belongs to the class 1A.

\smallskip
\noindent {\bf 2A case:} In this case, $\VFnX{2A}\simeq U_{6A}$.
Then $g=\tau_{e^1}\tau_{e^2}$ has order $2$ or $6$ and the group generated by
$\xi_u$ and $g$ is a cyclic group of order $6$ which is generated by
a 6A-element of $\M$.
Let $t$ be the unique involution in $\la \xi_u, g\ra$.
Then by \cite{ATLAS}, $t$ belongs to the 2A conjugacy class of $\M$
and
$C_\M(t)$ is isomorphic to a double cover of the Baby Monster $\B$.
Thus we have an exact sequence
\[
  1\longrightarrow \la t\ra \longrightarrow C_\M(t) \overset{\varphi_t} \longrightarrow  \B \longrightarrow 1.
\]
Since $\la t\ra$ and $\la \xi_u\ra$ are unique order 2 and order 3 subgroups in $\la \xi_u, g\ra$ and $t$ commutes with $\xi_u$, we have
\[
  N_\M( \la \xi_u, g\ra)= N_\M( \la \xi_u, t \ra) = C_\M(t) \cap N_\M(\la \xi_u\ra ).
\]
Set $G=N_\M( \la \xi_u, g\ra)$. Then
\[
\varphi_t( G)= N_{\B}(\la \varphi_t(\xi_u)\ra ) \quad \text{ and } \quad \psi_u(G)= C_{\Fi}(\psi_u(t)).
\]
Note that $\varphi_t(\xi_u)$ has order $3$ and $\psi_u(g)=\psi_u(t)$ has order $2$
since $(2,3)=1$.
By comparing the $3$-local subgroups of $\B$
and the $2$-local subgroups of $\Fi$ in \cite{ATLAS},
we have
\[
\varphi_t(G) \simeq \mathrm{S}_3\times \mathrm{Fi}_{22}\!:\!2
\quad \text{ and } \quad
\psi_u(G)=C_{\Fi}(\psi_u(g))\simeq (2\times 2.\mathrm{Fi}_{22})\!:\!2.
\]
Thus, $h=\psi_u(g)$ belongs to the conjugacy class 2A of $\Fi$ by \cite{ATLAS}.

\smallskip

\noindent {\bf 3A case:} In this case, $\VFnX{3A}\simeq M_{\mathcal{C}_4}$,
the ternary code VOA associated to the tetra code $\mathcal{C}_4$ and $\xi_u$,
$\tau_{e^1}$ and $\tau_{e^2}$ generate a subgroup of the shape $3^2\!:\!2$,
which has exactly $4$ distinct subgroups of order $3$.
Since $\VFnX{3A}\supset \W(\sfr{4}{5})^{\otimes 4}$ and each $\W(\sfr{4}{5})$
defines a non-trivial subgroup of order $3$, all order $3$ elements of
$\la \xi_u, \tau_{e^1}, \tau_{e^2}\ra$ belong to the conjugacy class 3A of $\M$.

Let $g=\tau_{e^1}\tau_{e^2}$. Then $\xi_u$ and $g$ generate a 3A-pure
elementary abelian $3$-subgroup of order $3^2$ in $\M$. The normalizer
$N_\M(\la \xi_u, g\ra )$ has the shape $(3^2\!:\!2\times \mathrm{O}^+_8(3)
).\mathrm{S}_4$ while the centralizer $C_\M(\la \xi_u,
g\ra)$ has the shape $3^2\times \mathrm{O}^+_8(3)$ (cf. \cite{Wi} and page 234 of \cite{ATLAS} ).
Thus,  $N_\M(\la \xi_u, g\ra )$ acts (by conjugation) on $\la
\xi_u, g\ra$ as $2\mathrm{S}_4(\simeq \mathrm{GL}_2(3))$ and $\mathrm{S}_4$ acts as
permutations of the $4$ distinct subgroups of order $3$ in $\la \xi_u, g\ra$.
Thus, $ N_\M(\la \xi_u\ra )\cap N_\M(\la \xi_u, g\ra)$ has the shape
\hbox{$(3^2\!:\!2\times \mathrm{O}^+_8(3)).\mathrm{S}_3$} and $\psi_u(
N_\M(\la \xi_u\ra )\cap N_\M(\la \xi_u,g\ra))$ has the shape
\hbox{$\mathrm{S}_3 \times \mathrm{O}^+_8(3).\mathrm{S}_3$.} Since
$\psi_u(\la \xi_u, g\ra) = \psi_u(\la g\ra)$,  we see that $\psi_u( N_\M(\la
\xi_u\ra )\cap N_\M(\la \xi_u, g\ra))$ normalizes the subgroup $\la \psi_u(g)\ra$
and hence $\psi_u( N_\M(\la \xi_u\ra )\cap N_\M(\la \xi_u, g\ra)) <  N_{\Fi}(\la
\psi_u(g)\ra)$.

By \cite{ATLAS}, page 207,  $\psi_u(N_\M(\la \xi_u\ra )\cap N_\M(\la \xi_u,
g\ra))$ is isomorphic to $N_{\Fi}(3A)$, where $N_{\Fi}(3A)$ denotes  the
normalizer of a cyclic subgroup generated by a 3A-element in $\Fi$, and it is a
maximal subgroup of $\Fi$. Thus, $N_{\Fi}(\la \psi_u(g)\ra) \simeq N_{\Fi}(3A)$
and $h=\psi_u(g)$ belongs to the conjugacy class 3A of $\Fi$. \qed

\begin{rem}
We note that in the 3A case, the group $N_\M(\la \xi_u, g\ra )$ acts on
$\VFnX{3A}\simeq M_{\mathcal{C}_4}$ as $(3^2:2)\mathrm{S}_4$, which is isomorphic to
$\aut(M_{\mathcal{C}_4})$ (cf. Remark \ref{rem:3.7}).
In fact, by the similar argument as in Remark \ref{32:2}, one can show without
using the property of the Monster that the whole group $\aut(M_{\mathcal{C}_4})$
can be extended to a subgroup of $\aut(V^\natural)$ since all irreducible modules
of $M_{\mathcal{C}_4}$ can be embedded into $V_{\sqrt{2}E_8}$ and are
$\aut(M_{\mathcal{C}_4})$-invariant.
\end{rem}


\appendix
\section{Appendix: Embedding of $\VFnX{3A}$ into $V^\natural$}\label{sec:A.4}

In this Appendix, we give  an embedding of $\VFnX{3A}$ into the moonshine VOA
which completes the proof of Proposition \ref{prop:5.29}. We achieve this by
providing an explicit embedding of $\VFnX{3A}\simeq M_{\mathcal{C}_4}$ into
$V_{\Lambda}^+\subset V^\natural$, where $M_{\mathcal{C}_4}$ refers to the
ternary code VOA associated to the tetra code $\mathcal{C}_4$ constructed in
\cite{KMY} (see Eq.~\eqref{eq:3.13+}). The main idea is essentially given  in
 \cite{GLSDC} and \cite{GL3}  (see also \cite{LM}).

\medskip

First, we consider some automorphisms of $V_{A_2}$ and $V_{E_8}$.
Set $a_i=E_{i,i}-E_{i+1,i+1}$, $x_i^+=E_{i,i+1}$ and $x_i^-=E_{i+1,i}$ for $i=1$, $2$,
where $E_{i,j}$ denotes a $3\times 3$ matrix whose $(i,j)$-entry is $1$ and others
are $0$.
Then $\{ a_i,\,x_i^\pm \mid  i=1,\,2\}$ is a set of Chevalley generators of
$\mathfrak{sl}_3(\C)$.
Let $\al_1$, $\al_2$ be simple roots of the root lattice $A_2$.
The weight one subspace of $(V_{A_2})_1$ of the lattice VOA
$V_{A_2}\simeq M_{\C A_2}(1)\otimes \C[A_2]$ forms a Lie algebra isomorphic
to $\mathfrak{sl}_3(\C)$ by the following correspondence:
\begin{equation}\label{eq:a.1}
  a_i\longmapsto {\alpha_i}_{(-1)}\vac,~~~
  x_i^\pm\longmapsto \pm e^{\pm \alpha_i}.
\end{equation}

The automorphism group of the lattice VOA $V_{A_2}$ is isomorphic to the
automorphism group of the Lie algebra $\mathfrak{sl}_3(\C)$, which is
isomorphic to $\mathrm{PSL}_3(\C)\rtimes \Z_2$.
Since $V_{A_2}$ is generated by its weight one subspace, $\mathrm{SL}_3(\C)$ acts
on $V_{A_2}$ via the adjoint map~\cite{FLM} under the identification above.
Namely, $P\in \mathrm{SL}_3(\C)$ acts by
$A\mapsto A^P=P^{-1}AP$ for $A\in \mathfrak{sl}_3(\C)\simeq (V_{A_2})_1$.

Let $\zeta=e^{2\pi i/3}$ be a cubic root of unity and consider the elements
\begin{equation}\label{tauands}
  \tau:=\left(
  \begin{array}{ccc}
    0 & 1 & 0
    \\
    0 & 0 & 1
    \\
    1 & 0 & 0
  \end{array}\right) \quad \text{ and } \quad
  s: =\frac{1}{\sqrt{3}}
  \left(\begin{array}{ccc}
    \zeta & \zeta^2 & 1
    \\
    \zeta^2 & \zeta & 1
    \\
    1 & 1 & 1
  \end{array}\right)
  \end{equation}
in $\mathrm{SU}_3\subset  \mathrm{SL}_3(\C)$.
Then
\begin{equation}\label{conj}
 r:= s^{-1} \tau  s =
  \left(\begin{array}{ccc}
    \zeta & 0 & 0
    \\
    0 & \zeta^2 & 0
    \\
    0 & 0 & 1
  \end{array}\right)_.
\end{equation}

Let $A_2^*$ be the dual lattice of $A_2$. We can describe the action of $r$ and
$\tau$ on $V_{A_2}$ and the $V_{A_2}$-module $V_{A_2^*}$ using this
identification more explicitly. Let $\delta =\al_1+\al_2$ be the half sum of positive
roots, i.e., the Weyl vector. Then $r(u\otimes e^\al) = \zeta^{\la \delta ,\al\ra } u
\otimes e^\al$ for all $u\in M(1)$ and $\al\in A_2^*$. Since $\tau$ normalizes
the canonical Cartan subalgebra of $\mathfrak{sl}_3(\C)$,  it induces an order
$3$ automorphism $\bar \tau$ of the root lattice $A_2$ explicitly given by
$$
  \ol{\tau}:~~ \al_1 \longmapsto \al_2
   \longmapsto -(\al_1+\al_2) \longmapsto \alpha_1.
$$

\medskip

Now fix a sublattice of type $A_2^4=A_2\perp A_2\perp A_2\perp A_2$ in $E_8$.
Then $E_8/A_2^4\subset (A_2^*/A_2)^4\simeq (\Z/3\Z)^4$ can be identified with the
tetra code ${\cal C}_4$.
The corresponding inclusion $V_{A_2}^{\otimes 4} \subset V_{E_8}$
induces an action of $\mathrm{SL}_3(\C)^4$ on $V_{E_8}$ by
automorphisms where the center $(\Z/3\Z)^4$ of $\mathrm{SL}_3(\C)^4$ acts via
its quotient $\widehat {\cal C}_4$, the dual group of ${\cal C}_4$.

Define
\[\begin{array}{rclcrcl}
\tilde{h}_1&=& 1\otimes \tau \otimes \tau \otimes \tau ,  &\qquad&
\tilde{h}_2&=& \tau \otimes \tau \otimes \tau^{-1} \otimes 1,\\
\rho_1&= &1 \otimes r \otimes r \otimes r , &\qquad&  \rho_2&=& r \otimes r \otimes r^{-1} \otimes 1, \\
\tilde{s}&=& s \otimes s\otimes s\otimes s &&
\end{array}\]
as automorphisms of $V_{E_8}$.

Let ${\cal G}$ be the subgroup generated by $\tilde{h}_1$ and $\tilde{h}_2$ and
$\mathcal{F}$ be the subgroup generated by $\rho_1$ and $\rho_2$.
Then $\mathcal{F}\simeq {\cal G}\simeq 3^2$. Moreover, by \eqref{conj},
we have
\begin{equation}\label{Hconj}
  \mathcal{F} = \tilde{s}^{-1}  {\cal G} \tilde{s}.
\end{equation}
\begin{rem}\label{rem:a1}
  We shall note that $\{x\in E_8\mid \la x, (0, \delta, \delta , \delta)\ra \in 3\Z\}
  \simeq E_6\perp A_2$ and  $\{x\in E_8\mid \la x, (\delta, \delta , -\delta, 0)\ra \in 3\Z\}
  \simeq E_6\perp A_2$. Moreover, we have
  \begin{equation}\label{eq:a5}
  K:=\{x\in E_8\mid \la x, (0, \delta, \delta , \delta)\ra \in 3\Z\text{ and } \la x, ( \delta, \delta ,-\delta, 0)\ra \in 3\Z \}\simeq A_2^4.
  \end{equation}
  Thus, we have $V_{E_8}^{\la \rho_1\ra}\simeq V_{E_8}^{\la \rho_2\ra} \simeq V_{E_6}\otimes V_{A_2}$ and $V_{E_8}^{\la \rho_1, \rho_2\ra} \simeq V_{A_2^4}$.
  However, we shall remark that the sublattice $K\simeq A_2^4$ obtained in \eqref{eq:a5} is not the same $A_2^4$ used
  to define~$h_i$ and~$\rho_i$, $i=1$,~$2$.
\end{rem}

Set
\[
  h_1=  \mathrm{id}\oplus  \bar\tau\oplus \bar \tau\oplus \bar\tau
  \quad \text{ and }\quad
  h_2=   \bar \tau\oplus \bar \tau\oplus \bar {\tau}^{-1} \oplus \mathrm{id}
\]
Then $h_1$ and $h_2$ can be considered as isometries of $E_8 \subset (A_2^*)^4$
induced by $\tilde{h}_1$ and $\tilde{h}_2$ and they also generate a subgroup
isomorphic to $3^2$.

\medskip

Now consider the following sublattices of an orthogonal sum $E_8\perp E_8$.
Let
\[
\begin{split}
  R & =\{ (x,x)\in E_8\perp E_8\mid x\in E_8\}, \\
  R^1 & =\{(x,h_1 x) \in E_8\perp E_8\mid x\in E_8\}, \\
  R^2 & =\{(x,h_2 x) \in E_8\perp E_8\mid x\in E_8\}.
\end{split}
\]
Then $R \simeq R^1\simeq R^2\simeq \sqrt{2}E_8$. Note also that
$R^1= (\mathrm{id}\oplus  h_1)R$ and $R^2 =(\mathrm{id}\oplus h_2)R$.

\smallskip

Let
\[
  e_R = \frac{1}{16} \w_R + \frac{1}{32} \sum_{\al\in R_4} e^\al
\]
be the Ising vector associated to $R$ defined by \eqref{eq:4.10}, where
$R_4=\{ \alpha \in R \mid \la \alpha,\alpha\ra =4\}$.

Now let $\widetilde{\mathcal{F}} = \{ \mathrm{id}\otimes \rho\mid  \rho\in \mathcal{F}\}$. Then
$\widetilde{\mathcal{F}}$ stabilizes the lattice VOA $V_R$ and by Remark~\ref{rem:a1}, the fixed point
sub-VOA $V_R^{\widetilde{\mathcal{F}}}\simeq V_{\sqrt{2}{A_2^4}}$. By the definition of
$\VFnX{3A}$ (see the 3A~case in Sec.~\ref{sec:3.3new}),  we can obtain a sub-VOA isomorphic to
$\VFnX{3A}$ in $V_R\simeq V_{\sqrt{2}E_8}$.  By Lemma~\ref{VF3A}, this $\VFnX{3A}$ is
generated by $e_R$,  $(\mathrm{id}\otimes \rho_1) e_R$ and
$(\mathrm{id} \otimes \rho_2) e_R$.

\medskip

Consider the automorphism  $\hat{\sigma}: = \ts\otimes \ts^{-1}$  of
$V_{E_8\perp E_8}\simeq V_{E_8} \otimes V_{E_8}$.
The following lemma is essentially proved in \cite{GL3} with some trivial modification.
\begin{lem}[Lemma 2.19 of \cite{GL3}]
  We have $\hat{\sigma}^{-1} e_R\in V_R$.
\end{lem}

\smallskip

Next, we will show that $\hat{\sigma}^{-1} e$ is in fact in $V_R^+$. For any
even lattice $L$, let $\theta:V_L\to V_L$ be the involution defined by
\begin{equation}\label{theta}
  \theta : \al_1(-n_1)\cdots \al_k(-n_k)\otimes e^\al
   \longmapsto (-1)^k \al_1(-n_1)\cdots \al_k(-n_k)\otimes e^{-\al}
\end{equation}
(cf. \cite{FLM,M1}).
Note that if $L$ is a root lattice of type $A_2$, by identifying
$(V_{A_2})_1$ with $\mathfrak{sl}_3(\C)$ as in \eqref{eq:a.1}
we have
\[
 A^\theta =- {}^t\! A, \qquad A\in \mathfrak{sl}_3(\C).
\]
If we extend $\theta$ to a map on $V_{A^*_2}$ then we have
$$ \mathrm{Aut}(\mathrm{SL}_3(\C)) \simeq \mathrm{SL}_3(\C)\rtimes \langle \theta \rangle$$
where $\theta$ acts by $X \mapsto {}^t\!X^{-1}$ on $\mathrm{SL}_3(\C)$.

\begin{lem}\label{plus}
  The Ising vectors $ \hat{\sigma}^{-1} e_R$,
  $(\mathrm{id}\otimes \tilde{h}_1) \hat{\sigma}^{-1} e_R$ and
  $(\mathrm{id} \otimes \tilde{h}_2) \hat{\sigma}^{-1}e_R$ are fixed by $\theta$.
\end{lem}

\pf
Let $s$ be defined as in \eqref{tauands}.
Then we have ${}^t s=s$ and $s^4=1$. Thus $s$ and $\theta$ generate
a dihedral group of order~$8$. One has
$\theta s\theta s^{-1}=\theta s^{-1} \theta s =s^2$ as a direct
calculation shows.

Now we consider the action of $\ts$ and $\theta$ on $V_{E_8}$.
This action is given by the diagonal embedding of $\langle s,\,\theta\rangle$
into $(\mathrm{Aut}(\mathrm{SL}_3(\C))^4$. We get
$\theta \ts\theta \ts^{-1}=\theta \ts^{-1} \theta \ts =\ts^2$.

Since $s^2$ is a permutation matrix, $s^2$  normalizes the
Cartan subalgebra $\C a_1+\C a_2$ of $\mathfrak{sl}_3(\C)$.
Moreover, $\ts^2$ normalizes the standard Cartan subalgebra of $(V_{E_8})_1$.
That means
\[
 \ts^2 (M_{\C E_8}(1)) = M_{\C E_8}(1).
\]
Thus $\ts^2$  induce an isometry $\mu:= \ol{\ts^2}$
of the root lattice $E_8$ such that
\begin{equation} \label{mu}
 \ts^2(e^\al) = \epsilon(\al) e^{\mu \al}
 \quad \text{ for }
 \al \in E_8,
\end{equation}
where $ \epsilon(\al)=\pm 1$.

Since $e^{(\alpha,\alpha)}={e^{(\alpha,0)}}_{(-1)}e^{(0,\alpha)}$ in
$V_{E_8\oplus E_8}\simeq V_{E_8}\tensor V_{E_8}$, one has
\[
  \theta\hat{\sigma}\theta\hat{\sigma}^{-1}(e^{(\al,\al)})
=\theta\hat{\sigma}\theta\hat{\sigma}^{-1}({e^{(\al,0)}}_{(-1)}e^{(0,\al)})
=(\theta\ts\theta\ts^{-1}({e^{(\al,0)}}))_{(-1)}
 (\theta\ts^{-1}\theta\ts({e^{(0,\al)}}))
\]
\[
 \qquad\qquad\qquad\qquad\qquad\qquad  = (\epsilon(\al) e^{(\mu \al, 0)})_{(-1)} (\epsilon(\al) e^{(0,\mu \al)})
  = e^{(\mu\al, \mu\al)}
\]
for any $\al \in E_8$.
Therefore, $\theta\hat{\sigma}\theta\hat{\sigma}^{-1}$ fixes $e_R$ and
so does $\hat{\sigma}\theta\hat{\sigma}^{-1}
=\theta(\theta\hat{\sigma}\theta\hat{\sigma}^{-1})$ since $\theta$ also fixes $e_R$.
Thus, $\hat{\sigma}^{-1}e_R$ is fixed by $\theta$.
Since $\mathrm{id}\otimes \tilde{h}_1$ and $\mathrm{id}\otimes \tilde{h}_2$
commute with $\theta$,  $(\mathrm{id}\otimes \tilde{h}_1)\hat{\sigma}^{-1}e_R$
and $(\mathrm{id}\otimes \tilde{h}_2) \hat{\sigma}^{-1} e_R$ are also fixed by
$\theta$.
\qed

\medskip

Using~(\ref{Hconj})  we have that $\hat{\sigma}^{-1} (\VFnX{3A})$ is generated by
$\{(\mathrm{id}\otimes \tilde{h}) \hat{\sigma}^{-1} e_R \mid \tilde{h} \in {\cal G} \}$ or
$\{ \hat{\sigma}^{-1} e_R,\, (\mathrm{id}\otimes \tilde{h}_1) \hat{\sigma}^{-1} e_R,\,
(\mathrm{id} \otimes \tilde{h}_2) \hat{\sigma}^{-1} e_R\}$.

Since $(\mathrm{id}\oplus h_1)R=R^1$ and $(\mathrm{id}\oplus h_2)R=R^2$,
we have $(\mathrm{id}\otimes \tilde{h}_1) \hat{\sigma}^{-1} e_R\in V_{R^1}^+$
and $(\mathrm{id}\otimes \tilde{h}_2) \hat{\sigma}^{-1} e_R \in V_{R^2}^+$, by
Lemma \ref{plus}.
Hence, we have $\hat{\sigma}^{-1}(\VFnX{3A}) < V_{R+R^1+R^2}^+$.
Therefore, it remains to show that $L=R+R^1+R^2$ can be embedded into
the Leech lattice~$\Lambda$.

\bigskip

First, we recall the ternary Golay code construction of the Leech lattice
$\Lambda$ \cite{cs}.

Let $S$ be an orthogonal sum of $12$ copies of $A_2$.
Then the discriminant group $S^*/S$ has a natural identification
with $\Z_3^{12}$.
The ternary Golay code ${\cal C}_{12}\subset \Z_3^{12}$ can be defined using
the tetracode ${\cal C}_4$ as the set
$${\cal C}_{12}=\{(c^0,-c^+,c^-)\mid c^0,c^+,c^-\in \Z_3^4,\
c^0+c^++c^-=-\hbox{$\sum_{i=1}^4c^0_i$}\cdot (1,1,1,1),\ c^+-c^-\in {\cal C}_4 \}.$$
For each codeword $x=(x_1,\dots, x_{12})\in {\cal C}_{12}$, let
$ \gamma_x = (\gamma_{x_1}, \dots, \gamma_{x_{12}})\in S^*$ be some vector
which modulo $S$ gives the codeword $x$. Then
\[
\mathcal{N} : = \bigcup_{{x \in \cal C}_{12}}  (\gamma_x +S)
\]
is isometric to the Niemeier lattice of type $A_2^{12}$.

Let $\delta= \al_1+\al_2$ be the half sum of positive roots of $A_2$ and let
$$\hat{\delta}: =(\delta,\delta, \delta,\delta, -\delta, -\delta, -\delta, -\delta,
\delta,\delta, \delta,\delta).$$
Then
\[
  \mathcal{N}^0=\{ \al\in \mathcal{N} \mid (\al, \hat{\delta})\in 3\Z\}
\]
is a sublattice of index~$3$ without roots.
Note that $(\al_1, 0,\dots,0)+\frac{1}3 \hat{\delta}$ has  norm $4$ and the
lattice $ \mathcal{N}^0 + \Z((\al_1, 0,\dots,0)+\frac{1}3 \hat{\delta})$ is even
unimodular without roots.
Hence, it is isometric to the Leech lattice $\Lambda$ (see Chapter 24 of \cite{cs}).

\medskip

Next, we construct some $\sqrt{2}E_8$ sublattices of $\mathcal{N}^0 < \Lambda$.
Note that
$$
  \{(0^4,c,c)\mid c\in {\cal C}_4\} < {\cal C}_{12}.
$$
Thus
$$
  \tilde{R}:=\mathrm{span}\{(0,\,\gamma_c+z,\,\gamma_c+z)\mid z\in A_2^4,\
  c\in {\cal C}_4\}<\mathcal{N}.
$$
Since $E:=\bigcup_{c\in {\cal C}_4} (\gamma_c+ A_2^4) \simeq E_8$
it follows that $\tilde{R} \simeq \sqrt{2} E_8$.

Let $E^1:=0\oplus E \oplus 0$ and $E^2:= 0\oplus 0 \oplus E$.
Then $(E^1, E^2)=0$ and $E^1 + E^2 <  \frac{1}3 \Lambda$.

For a codeword $x\in {\cal C}_{12}$,
let $h(x):= \bar \tau^{x_1} \oplus \cdots \oplus \bar \tau^{x_{12}}$
where $\bar \tau$ is the previously defined isometry of $A_2$.
Then $h(x)$ is an isometry of $\mathcal{N}$  and $\Lambda$ \cite{cs}.

Consider now the codewords
$$
  d^1=(0,-1,-1,-1,~0,0,0,0,~0,1,1,1)\quad \hbox{and}\quad
  d^2=(1,1,0,-1,~0,0,0,0,~1,1,-1,0)
$$
of ${\cal C}_{12}$.
Define $\hat h_1:=h(d^1)$  and  $\hat h_2:=h(d^2)$.
Note that $\hat h_1$ and $\hat h_2$  act on $E^1+ E^2$ as
$\mathrm{id} \oplus h_1$ and $\mathrm{id} \oplus  h_2$, where
${h}_1= \mathrm{id} \oplus \bar \tau \oplus \bar \tau \oplus \bar \tau $,
and
${h}_2= \bar \tau \oplus \bar \tau \oplus \bar \tau^{-1} \oplus \mathrm{id}$
in $O(E_8)$ as previously defined.
Then
\[
\begin{split}
\tilde{R^1}= \hat h_1(\tilde{R}) = \{ (0,\al, \tilde h_1 \al)\mid \al \in E^1\},\\
\tilde{R^2}= \hat h_2(\tilde{R}) = \{ (0,\al, \tilde h_2 \al)\mid \al \in E^1\}
\end{split}
\]
are contained in $\mathcal{N}^0< \Lambda$.
It is also clear that $\tilde L=\tilde{R}+\tilde{R^1}+\tilde{R^2}$ is
isometric to $L=R+R^1+R^2$.

\begin{rem}
It is known that a $3$-transposition group generated by three involutions $t_1$,
$t_2$, $t_3$ such that $t_3\notin\la t_1, t_2\ra$ and any two of them generate
$\mathrm{S}_3$ is either isomorphic to $\mathrm{S}_4$, $3^{1+2}{:}2$ or $3^2{:}2$ 
(see Lemma~2.5 of \cite{CH1}). 
By Remark~\ref{32:2}, the $\tau$-involutions associated to 
$\{ \hat{\sigma}^{-1} e_R,\, (\mathrm{id}\otimes \tilde{h}_1) \hat{\sigma}^{-1}
e_R,\, (\mathrm{id} \otimes \tilde{h}_2) \hat{\sigma}^{-1} e_R\}$ generate a
group of the shape $3^2{:}2$ in $\aut(V_\Lambda^+)$ and in $\aut(V^\natural)$. 
Finally, we note that  the centralizer of a 3C-element in $\mathrm{Co}_0$ has 
the shape $3^{1+4}{:} \mathrm{Sp}_4(3) \times 2$.  
It has a natural subgroup $3^{1+2}{:}2$  which is  generated by involutions with 
trace $-8$ on the Leech lattice (see (10.35.3) of \cite{Gr} or \cite{ATLAS}).  
Recall that an involution of trace $-8$ on $\Lambda$ has the fixed sublattice 
isometric to $\sqrt{2}E_8$ (cf.~Theorem (10.15) of \cite{Gr} or Chapter~10 of 
\cite{cs}). 
Thus, by formula \eqref{eq:4.10}, one can construct Ising vectors in 
$V_\Lambda^+ < V^\natural $ such that the corresponding Miyamoto involutions 
generate a group  of the shape $3^{1+2}{:}2$ in the Monster and the product of 
any two of them is in the conjugacy class 3A. 
In this case, the sub-VOA generated by the corresponding Ising vectors will not 
be isomorphic to $M_{\mathcal{C}_4}$ but the authors do not know the exact 
structure of such a VOA. 
\end{rem}

\small

\end{document}